\documentclass[nonblindrev, opre]{informs3}
\OneAndAHalfSpacedXI 

\usepackage{natbib}
\bibpunct[, ]{(}{)}{,}{a}{}{,}%
\def\bibfont{\small}%

\usepackage{graphicx}
\usepackage{xcolor}
\usepackage{soul}
\usepackage{multirow}
\usepackage{algorithm,algpseudocode}
\usepackage{bm}
\usepackage{subfig}

\TheoremsNumberedThrough     
\EquationsNumberedThrough    

\usepackage{fancyhdr}
\usepackage{xspace}
\usepackage{mathtools}

\usepackage{tikz}
\usetikzlibrary{calc}
\usetikzlibrary{decorations.pathreplacing}
\usetikzlibrary{arrows,automata}
\usetikzlibrary{positioning}
\usetikzlibrary{shapes}

\usepackage{ctable}
\newcommand{\Path}{\rho}
\newcommand{\Graph}{\mathcal{G}}
\newcommand{\Arc}{\mathcal{A}}
\newcommand{\Edge}{\mathcal{V}}

\DeclarePairedDelimiter\abs{\lvert}{\rvert}%
\DeclarePairedDelimiter\norm{\lVert}{\rVert}%

\makeatletter
\let\oldabs\abs
\def\abs{\@ifstar{\oldabs}{\oldabs*}}
\let\oldnorm\norm
\def\norm{\@ifstar{\oldnorm}{\oldnorm*}}
\makeatother

\definecolor{light-gray}{gray}{0.95}

\begin{document}

\RUNAUTHOR{Schrotenboer et al.}
\RUNTITLE{Reliable Reserve-Crew Scheduling for Airlines}
\TITLE{{Reliable} Reserve-Crew Scheduling for Airlines}
\ARTICLEAUTHORS{%
	\AUTHOR{Albert H. Schrotenboer}
	\AFF{Operations, Planning, Accounting \& Control Group, School of Industrial Engineering, Eindhoven University of Technology\\ \EMAIL{a.h.schrotenboer@tue.nl}}
	\AUTHOR{Rob Wenneker}
	\AFF{Space, Transport \& Logistics, CGI Netherlands, Rotterdam}
	\AUTHOR{Evrim Ursavas, Stuart X. Zhu}
	\AFF{Department of Operations, Faculty of Economics and Business, University of Groningen}

} 

\ABSTRACT{We study the practical setting in which regular- and reserve-crew schedules are dynamically maintained up to the day of executing the schedule. At each day preceding the execution of the schedule, disruptions occur due to sudden unavailability of personnel, making the planned regular and reserve-crew schedules infeasible for its execution day. This paper studies the fundamental question on how to repair the schedules' infeasibility in the days preceding the execution, taking into account labor regulations. We propose a reliable repair strategy that maintains flexibility in order to cope with additional future disruptions. The flexibility in reserve-crew usage is explicitly considered through evaluating the expected shortfall of the reserve-crew schedule based on a Markov chain formulation. The core of our approach relies on iteratively solving a set-covering formulation, which we call the reliable Crew Recovery Problem, which encapsulates this flexibility notion for reserve crew usage. A tailored branch-and-price algorithm is developed for solving the reliable Crew Recovery Problem to optimality. The corresponding pricing problem is efficiently solved by a newly developed pulse algorithm. Based on actual data from a medium-sized hub-and-spoke airline, we show that embracing our approach leads to fewer flight cancellations and fewer last-minute alterations, compared to repairing disrupted schedules without considering our robust measure.}
\KEYWORDS{Airline Operations, Disruption Management, Branch-and-Price Algorithm, Reserve Crew, Scheduling}

\maketitle
\section{Introduction}
Due to a significant growth in air traffic, airports are becoming more and more congested \citep{jacquillat2015}. Consequently, efficient and effective disruption management is becoming crucial for hub-and-spoke airlines in order to stay competitive \citep{iata2011vision}. The airline's capability to deal with inevitable disruptions such as crew absenteeism and sudden aircraft unavailability depends predominantly on the flexibility of an airline's (reserve) crew schedule. The literature on how to create a reliable crew schedule during initial schedule creation \citep[see, e.g.,][]{klabjan2002airline, yen2006stochastic, cacchiani2016optimal, wei2018modeling}, and how to recover from perturbed schedules on the day of execution \citep[see, e.g.,][]{rosenberger2003rerouting,petersen2012optimization, maher2015, ruther2016integrated}, is developing. However, the effect of using reserve crew members on the robustness of airline crew schedules, and how this affects the airline's capability to deal with future disruptions in the days before schedule execution, is still unknown \citep{wen2021airline}. We will study these interactions where we incorporate a robust measure for the use of specific reserve-crew schedules in an airline crew recovery problem. In this way, a more flexible schedule is provided which leads to fewer flight delays, flight cancellations, and changes in the crew schedule. 

In practice, crew schedules are made available to the crew several weeks before the flights are executed and are continuously updated to ensure that all operations continue as planned. In the period between schedule publishing and flight execution, known as the \textit{tracking period}, unexpected events may occur, which affects an airline's capability to deal with future disruptions. To deal with such disruptions, airlines typically exploit a reserve crew schedule next to their actual operations. However, using a reserve crew-member might have severe implications for the ability to recover from future disruptions in the days until the execution of the schedule. In other words, it is of crucial importance that recovering from disrupted schedules is done in such a way that the flexibility (in reserve-crew usage) is preserved for recovering from future disruptions. Very little research has been done till now on how to maintain a reliable schedule in the tracking period that can cope with additional disruptions in the hours or days following current disruptions. 

In our study, we consider a set of so-called crew-induced disruptions occurring during the tracking period. A crew-induced disruption is defined as the inability of a crew member to follow their schedule consisting of flight legs (origin to destination flights) and reserve shifts. Those may be due to reasons such as crew taking days off, sickness, as well as crew going over their maximum duty time. In short, we consider reasons that a crew member becomes unavailable to fly on a given day in the schedule. Through re-timing flights, swapping crew, and using reserve crew members, the goal is to return to a feasible, reliable, and cost-efficient schedule, while preventing unnecessary altering of published schedules and respecting labor regulations. Robustness is explicitly considered through the penalization of unfavorable pairing characteristics and through evaluating the effect of using reserve crew based on the underlying reserve crew schedule. Basically, our approach modifies the cost of a crew-pairing integer program to account for the expected shortfall of personnel on the day of flight execution for selecting a specific reserve-crew schedule. This is contrary to the traditional crew recovery literature, where reserve crew costs are usually generalized to a single cost parameter which is given exogenously rather than depending on the actual reserve-crew schedule being used.

Our approach is also different from the well-known robustness approaches in the literature on robustness optimization {}{\cite[see, e.g.,][]{Soyster73, Ben98, Bertsimas04}}.
Under the assumption that the uncertain data resides in the so-called uncertainty set, {}{these approaches try to find a feasible solution for all possible scenarios or deviations of uncertain parameters in their uncertainty set} and aim at minimizing the objective function in the worst case. Interested readers can refer to \cite{Bertsimas11} for a comprehensive review. In our study, we modify the cost of a standard integer program to balance between standard air-crew recovery cost and the expected shortfall of personnel of selecting a particular reserve-crew schedule. Thus, our approach aims to create reliable reserve-crew schedules to control the risk of flight cancellations or last minute flight alterations but does not employ standard robust optimization techniques.

A Markov chain framework is presented to evaluate reserve crew schedules underlying the flight schedule for a single day. We evaluate reserve-crew schedules based on the Earliest Finisher First (EFF) recovery policy, in which reserve-crew members are used based on their remaining time until the end of their reserve duty. We prove that under some conditions this recovery policy is optimal. We test the effect of the developed robust measure for reserve-crew usage by repairing disrupted schedules dynamically during the tracking period. The associated optimization problem, which needs to be solved multiple times during the tracking period, is called the reliable Crew Recovery Problem (RCRP). 

We present a set-covering model for the RCRP, and solve it using a branch-and-price algorithm. The associated pricing problem, which can be reduced to the Resource Constrained Shortest Path Problem \cite[see e.g.][]{irnich2005}, is solved using a tailored pulse algorithm \cite[see, e.g.,][]{lozano2015, schrotenboer2018}. We show through a case study from a medium-sized hub-and-spoke carrier in the Netherlands that our reliable approach outperforms traditional crew recovery methods. In particular, our approach results in a lower risk of having an insufficient number of reserve crew members, which leads to fewer flight cancellations and alterations.

In the following, we will highlight this paper's contributions by reviewing the relevant literature on airline recovery operations and airline (reserve) crew scheduling. We first review the crew pairing problem \cite[see, e.g.,][]{barnhart2003} and then discuss the crew recovery problem \cite[see, e.g.,][]{maher2015} afterward.



The crew pairing problem forms the basis of the initial schedule creation. It consists of generating minimum-cost anonymous multiple-day work schedules that satisfy legal and contractual obligations \citep{barnhart2003}. Several papers have solved the deterministic crew pairing problem to optimality using branch-and-price \citep[e.g.][]{desaulniers1997, gamache1999}. \cite{Zeighami19} propose a model that integrates the crew pairing and crew assignment problems simultaneously for pilots and copilots. They develop a method that combines {Benders} decomposition and column generation to solve the model. \cite{Quesnel20} consider crew rostering problem with crew preferences. They solve the problem by using a column generation algorithm in which new pairings are generated by solving subproblems consisting of constrained shortest path problems. However, in the creation of those pairings, robustness was not considered and the resulting schedules were extremely fragile to the effects of disruptions during future operations. Since then, reliable approaches to this crew pairing problem are presented. For example, \citet{shebalov2006} aim to maximize the potential number of so-called move-up crews, who are able to move to subsequent later flights if necessary. This allows for mitigating further delays by swapping crew members from different flights. \citet{schaefer2005} discuss a model in which unfavorable characteristics of a crew pairing, such as a crew's duty time being close to the maximum allowed duty time, are discouraged using a penalty function. The robust measure for flexibility of reserve-crew members, as is included in the RCRP, follows this approach by providing a sophisticated penalty function (based on a Markov-chain formulation) to increase our flexibility to cope with future disruptions. Recently, \cite{antunes2019robust} present a robust pairing model for airline crew scheduling, but as opposed to our approach, they do not consider assigning crew pairings to crew members and therefore do not consider the creation of reserve-crew schedules.

When disruptions are impossible to be covered by the regular crew, reserve crew on call needs to be utilized. One of the first authors to consider the reserve crew pairing problem were \citet{dillon1999}. They consider reserve crew as a separate entity and use a set partitioning formulation to create multi-day reserve pairings which cover the expected demand for reserve personnel. A similar problem is introduced by \citet{sohoni2006}, who consider the reserve cockpit crew scheduling from the perspective of an American airline in which the authors regard open trips resulting from bidline-invoked conflicts as the primary reason for reserve demand. More recently, \citet{bayliss2016} discusses the airline reserve crew scheduling problem for hub- and spoke airlines. He proposes several probabilistic models to evaluate reserve crew schedules which can be used to mitigate crew delay- and absence risk. Different from previous studies, we include the evaluation of reserve schedules in a crew recovery setting, where future schedule flexibility is explicitly considered. 

Crew recovery is first introduced by \citet{wei1997}. In their seminal paper, the authors consider the repair of disrupted crew schedules using a set partitioning formulation of possible repair actions in a branch-and-price framework. Whereas \cite{wei1997} do not allow for cancelling flights, \citet{lettovsky2000} do allow this and thereby extend the seminal work of \cite{wei1997}.  \citet{stojkovic1998} are the first to incorporate airline crew recovery in a personalized schedule setting. Their objective is to minimize the recovery costs with the additional consideration that the adjusted monthly assignment should remain as close to the original schedule as possible, as often schedule changes affects crew happiness.

The early literature on airline recovery provides sequential solution methods. First, aircraft recovery is performed for a given disruption. Subsequently, for the resulting aircraft schedule, crew recovery takes place. This approach is shown by \citet{papadakos2009} to generally be sub-optimal and not even necessarily feasible. In the past decades, numerous authors proposed integrated recovery models integrating aircraft- and crew recovery. \citet{stojkovic2001} improve on their crew recovery model by including fixed time windows for which an aircraft may be delayed and provide an integrated pilot and aircraft recovery problem. \citet{maher2015} provides a new approach to the integrated one-day aircraft and crew recovery problem using column- and row generation. He argues that, with this improved solution method, large integrated problems can be solved efficiently and quickly. \citet{Wen20} study a crew pairing problem with flight flying time variability.
They develop a reliable scheduling by encouraging deviation-affected-free flights and discouraging deviation-affected flights. A customized column generation based solution algorithm is proposed.
\citet{Wang20} formulate a two-stage stochastic program with integer recourse to 
jointly optimize scheduling interventions and ground-holding operations across airports networks under operating uncertainty. They develop an original decomposition algorithm to solve it.
\citet{Xu21} investigate a robust scheduling problem of schedule design, fleet assignment, and aircraft routing by incorporating propagated delays and flight re-timing decisions. A column generation procedure as well as a sequential variable neighborhood search heuristic are designed for solution. In our paper, we consider an integrated approach by allowing for the re-timing of flights, which increases the possibility to recover from disruptions significantly, as is shown by \citet{mercier2007}. 

Recently, the tracking period has been the motivation for \cite{ruther2016integrated} to study an integrated airline scheduling problem. To the best of the authors' knowledge, this is the only work focusing on the tracking period. However, those authors do not consider the actual creation of reserve-crew schedules but instead integrate other aspects of airline operations into their optimization problems. We refer {interested readers} on the integration of other airline operations into the crew pairing or crew recovery problem to the excellent review by \cite{kasirzadeh2017airline}, refer to \cite{woo2021scenario} for airline rescheduling problems from the perspective of a single airport, and refer to \cite{hu2016integrated} for the recovery of aircraft and passengers after airline operations.

Summarizing, the contributions of this paper are twofold. First, we explicitly model the effects of using reserve crew members on the robustness of the crew schedule. Based on a Markov chain formulation of the underlying reserve crew schedule, we are able to provide a penalty function ensuring robustness with regards to reserve crew usage. Our second contribution is to incorporate this notion of reliable reserve crew usage within an airline crew recovery model, which, next to recovering the disrupted schedule in a cost-efficient manner, has the secondary objective to maintain a schedule that is capable of absorbing further disruptions. Here we show that using our approach instead of traditional crew recovery methods is more reliable.
In short, we use Table \ref{tab:review} to highlight our contributions to the existing literature.

The remainder of this paper is structured as follows. In Section \ref{sec:propmodel}, we formally describe the RCRP. In Section \ref{sec:coststructure}, we introduce the cost structures for the crew pairings and the reserve-crew schedules. For the latter, we discuss a framework to evaluate reserve-crew schedules and their effect on the robustness of the overall crew schedule. The branch-and-price approach, including the pulse algorithm for solving the pricing problems, is discussed in Section \ref{sec:pricingproblem}. Section \ref{sec:experimentaldesign} will present detailed results on actual data of a medium-sized hub and spoke carrier. We conclude our work, and provide avenues for further research, in Section \ref{sec:conclusion}.

\begin{table}[htbp]
\begin{center}
{\footnotesize
\colorbox{red!20}{
\begin{tabular}
{l|c|c|c|c|c|c}\hline
Literature &  FA & AR & CS & GO & Model & Algorithm\\\hline
\citet{schaefer2005} &  &   & $\surd$ &               & MC & DA \\
\citet{shebalov2006} &  &   & $\surd$ &               & MIP& RO \\
\citet{papadakos2009} & $\surd$  &  $\surd$ & $\surd$ &               & MIP& SO \\
\citet{maher2015}  & $\surd$ &   & $\surd$ &               & MIP & CRG \\
\citet{bayliss2016} & $\surd$ &   & $\surd$ &               & MIP & SO\\
\cite{hu2016integrated}& $\surd$ &   & $\surd$ &               & IP & GRASP  \\
\cite{ruther2016integrated}&  & $\surd$  & $\surd$ &               & MIP & B\&P \\
\cite{antunes2019robust}  &  &   & $\surd$ &               & MIP & RO \\ 
\cite{Zeighami19}  &  &   & $\surd$ &               & MIP & BD, CG \\ 
 \cite{Quesnel20} &  &   & $\surd$ &               & IP & CG \\ 
\citet{Wang20} &  $\surd$ &  &  & $\surd$               & SP & BD \\
\citet{Wen20} &  &   & $\surd$ &               & MIP & CG \\
\citet{Xu21} &  $\surd$ & $\surd$  &  &               & MIP & CG \\
\cite{woo2021scenario} &  $\surd$ &  &$\surd$   &               & SP & GH \\
Our paper   & $\surd$  &  & $\surd$ &  & MIP, MC & B\&P \\ \hline
 \multicolumn{7}{l}{ FA: Fleet Assignment,  AR: Aircraft Routing, CS: Crew Scheduling, GO: Ground Operations } \\
\multicolumn{7}{l}{IP: Integer Programming, MIP: Mixed Integer Programming, SP: Stochastic Programmming, MC: Markov Chain} \\ 
\multicolumn{7}{l}{B\&B: Branch-and-Bound, B\&P: Branch-and-Price, 
BD: Benders Decomposition, CG: Column Generation} \\  
\multicolumn{7}{l}{ CRG: Column-and-Row Generation, DA: Deterministic Approximation, RO: Robust Optimization } \\ 
\multicolumn{7}{l}{GRASP: Greedy Randomized Adaptive Search Procedure, GH: Greedy Heuristic, S: Simulation Optimization}\\ 
   \hline
\end{tabular}}
}
\end{center}
\caption{\small{Related literature} }
\label{tab:review}
\end{table}

\section{System description} \label{sec:propmodel}
In this section, we formally describe the system that embeds the reliable Crew Recovery Problem (RCRP). First, we provide a general problem description of the RCRP, and second, we provide a Mixed Integer Programming (MIP) formulation. However, we start this section with an illustrative example of our system and the RCRP, and sketch the practical context of airline scheduling operations.

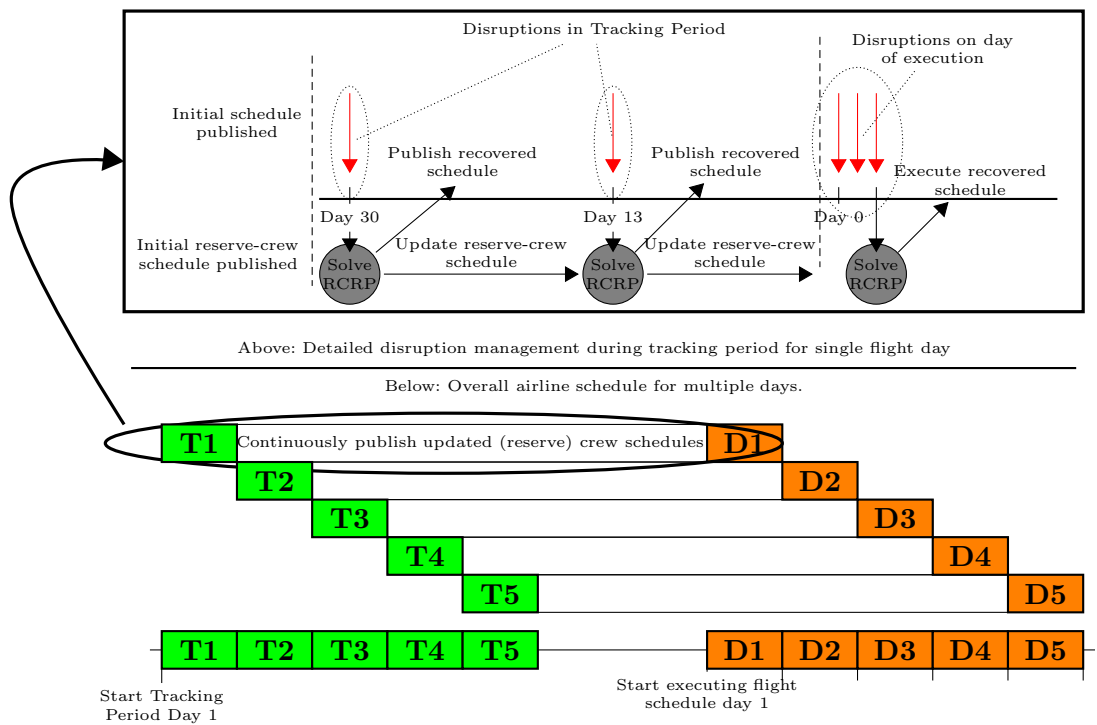
\begin{figure}[!htb]
\centering
\begin{tikzpicture}
\tiny
\node (v1) at (-5,0) {};
\node (v2) at (5,0) {};
\draw [thick] (v1) edge (v2);
\node (v3) at (-4.5,1.5) {};
\node (v4) at (-1,1.5) {};
\node (v5) at (-1,0) {};
\node (v10) at (-4.5,0.25) {};
\node (v11) at (-4.5,-0.25) {Day 30};
\node (v12) at (-1,0.25) {};
\node (v13) at (-1,-0.25) {Day 13};
\draw  (v10) edge (v11);
\draw  (v12) edge (v13);
\draw [-triangle 60, red] (v3) edge (v10);
\draw [-triangle 60, red] (v4) edge (v12);
\node (v20) at (-5,2) {};
\node (v21) at (-5,-1.25) {};
\node (v18) at (1.75,2.25) {};
\node (v19) at (1.75,-1) {};
\draw [densely dashed] (v18) edge (v19);
\draw [densely dashed] (v20) edge (v21);
\node (v30) at (-1.25,2.25) {Disruptions in Tracking Period};
\node (v22) at (2,0.25) {};
\node (v23) at (2,-0.25) {Day 0};
\draw  (v22) edge (v23);
\node[align = center] (v34) at (3.25,2) {Disruptions on day \\ of execution};
\node (v24) at (2,1.5) {};
\node (v25) at (2.25,1.5) {};
\node (v27) at (2.5,1.5) {};
\node (v28) at (2.5,0.25) {};
\node (v26) at (2.25,0.25) {};
\node at (2.25,1.5) {};
\node at (2.25,0.25) {};
\draw [-triangle 60, red] (v24) edge (v22);
\draw [-triangle 60, red] (v25) edge (v26);
\draw [-triangle 60, red] (v27) edge (v28);
\node[align=center] at (-6,1) {Initial schedule \\ published};

\node[align = center] at (-6.25,-0.75) {Initial reserve-crew \\ schedule published};
\draw [very thick] (-7.5,2.5) rectangle (5.25,-1.5);
\draw [densely dotted] (2.25,0.75) node (v35) {} ellipse (0.6 and 1);

\draw [densely dotted] (v35) edge (v34);
\draw [fill = gray] (-1,-1) ellipse (0.4 and 0.4);
\draw[fill = gray ](2.5,-1) ellipse (0.4 and 0.4);
\draw[fill = gray]  (-4.5,-1) ellipse (0.4 and 0.4);
\node [align = center] (v36) at (-4.5,-1) {Solve \\ RCRP};

\node (v38) at (-3,0.25) {};
\node[align = center] (v37) at (-1,-1) {Solve \\ RCRP};
\draw [-triangle 60, pos = 0.6] (v36) edge  node[align = center,auto, sloped,pos=0.5] {Update reserve-crew \\ schedule} (v37);
\draw [-triangle 60] (v36) edge (v38);

\node[align = center] at (-3,0.5) {Publish recovered \\ schedule};
\draw [-triangle 60] (v11) edge (v36);
\draw [-triangle 60] (v13) edge (v37);
\node (v39) at (1.75,-1) {};
\draw [-triangle 60] (v37) edge node[align = center,auto, sloped,pos=0.5] {Update reserve-crew \\ schedule} (v39);
\node[align = center] (v40) at (0.5,0.5) {Publish recovered \\ schedule};
\draw [-triangle 60] (v37) edge (v40);
\node[align=center] (v6) at (2.5,-1) {Solve \\ RCRP};

\node[align = center] (v7) at (3.75,0.25) {Execute recovered  \\ schedule};
\draw [-triangle 60] (v6) edge (v7);

\draw [-triangle 60] (v28) edge (v6);
\draw [densely dotted] (-4.5,0.75) node (v8) {} ellipse (0.25 and 0.75);
\draw [densely dotted] (-1,0.75) node (v9) {} ellipse (0.25 and 0.75);
\draw [densely dotted] (v8) edge (v30);
\draw [densely dotted] (v9) edge (v30);

\draw  (-7,-3) node (v63) {} rectangle (1.25,-3.5) node (v14) {};
\draw  (-6,-3.5) node (v64) {} rectangle (2.25,-4) node (v15) {};
\draw  (-5,-4) node (v65) {} rectangle (3.25,-4.5) node (v16) {};
\draw  (-4,-4.5) node (v66) {} rectangle (4.25,-5) node (v17) {};
\draw  (-3,-5) node (v67) {} rectangle (5.25,-5.5) node (v29) {};
\draw [fill = orange, thick] (0.25,-3) rectangle (1.25,-3.5);
\draw [fill = orange, thick] (v14) rectangle (v15);
\draw [fill = orange, thick] (v15) rectangle (v16);
\draw [fill = orange, thick] (v16) rectangle (4.25,-5);
\draw [fill = orange, thick] (v17) rectangle (v29);
\node (v31) at (-7.25,-6) {};
\node (v32) at (5.5,-6) {};
\draw  (v31) edge (v32);
\node (v33) at (-7,-5.75) {};
\node[align=center] (v41) at (-7,-6.75) {Start Tracking \\ Period Day 1};
\node (v42) at (-6,-5.75) {};
\node (v43) at (-6,-6.25) {};
\node (v44) at (-5,-5.75) {};
\node (v45) at (-5,-6.25) {};
\node (v46) at (-4,-5.75) {};
\node (v47) at (-4,-6.25) {};
\node (v48) at (-3,-5.75) {};
\node (v49) at (-3,-6.25) {};
\node (v50) at (0.25,-5.75) {};
\node[align = center] (v51) at (0.25,-6.6) {Start executing flight \\ schedule day 1};
\node (v52) at (1.25,-5.75) {};
\node (v53) at (1.25,-6.5) {};
\node (v54) at (2.25,-5.75) {};
\node (v55) at (2.25,-6.5) { };
\node (v57) at (3.25,-5.75) {};
\node (v56) at (3.25,-6.5) { };
\node (v58) at (4.25,-5.75) {};
\node (v59) at (4.25,-6.5) { };
\node (v60) at (5.25,-5.75) {};
\node (v61) at (5.25,-6.5) { };
\draw  (v33) edge (v41);
\draw  (v42) edge (v43);
\draw  (v44) edge (v45);
\draw  (v46) edge (v47);
\draw  (v48) edge (v49);
\draw  (v50) edge (v51);
\draw  (v52) edge (v53);
\draw  (v54) edge (v55);
\draw  (v56) edge (v57);
\draw  (v58) edge (v59);
\draw  (v60) edge (v61);

\node (v62) at (-2.9,-3.25) {Continuously publish updated (reserve) crew schedules};
\node at (0.75,-3.25) {\normalsize\textbf{D1}};
\node at (1.75,-3.75) {\normalsize\textbf{D2}};
\node at (2.75,-4.25) {\normalsize\textbf{D3}};
\node at (3.75,-4.75) {\normalsize\textbf{D4}};
\node at (4.75,-5.25) {\normalsize\textbf{D5}};
\draw[very thick]  (-3.25,-3.25) ellipse (4.5 and 0.4);

\draw[very thick, -triangle 60]  plot[smooth, tension=.7] coordinates {(-7.5,-3) (-9,0) (-7.5,0.5)};
\draw [fill = orange, thick] (0.25,-5.75) rectangle (1.25,-6.25);
\draw [fill = orange, thick] (v52) rectangle (2.25,-6.25);
\draw [fill = orange, thick] (v54) rectangle (3.25,-6.25);
\draw [fill = orange, thick] (v57) rectangle (4.25,-6.25);
\draw [fill = orange, thick] (v58) rectangle (5.25,-6.25);
\node at (0.75,-6) {\normalsize\textbf{D1}};
\node at (1.75,-6) {\normalsize\textbf{D2}};
\node at (2.75,-6) {\normalsize\textbf{D3}};
\node at (3.75,-6) {\normalsize\textbf{D4}};
\node at (4.75,-6) {\normalsize\textbf{D5}};
\node at (2.5,-7.5) {};
\draw [fill = green, thick] (v33) rectangle (v43);
\draw [fill = green, thick] (-6,-5.75) rectangle (v45);
\draw [fill = green, thick] (v44) rectangle (v47);
\draw [fill = green, thick] (-4,-5.75) rectangle (-3,-6.25);
\draw [fill = green, thick] (v48) rectangle (-2,-6.25);
\node at (-6.5,-6) {\normalsize\textbf{T1}};
\node at (-5.5,-6) {\normalsize\textbf{T2}};
\node at (-4.5,-6) {\normalsize\textbf{T3}};
\node at (-3.5,-6) {\normalsize\textbf{T4}};
\node at (-2.5,-6) {\normalsize\textbf{T5}};

\draw [fill = green, thick] (v63) rectangle (v64);
\draw [fill = green, thick] (v64) rectangle (-5,-4);
\draw [fill = green, thick] (v65) rectangle (v66);
\draw [fill = green, thick] (v66) rectangle (-3,-5);
\draw [fill = green, thick] (v67) rectangle (-2,-5.5);

\node[black] at (-6.5,-3.25) {\normalsize\textbf{T1}};
\node at (-5.5,-3.75) {\normalsize\textbf{T2}};
\node at (-4.5,-4.25) {\normalsize\textbf{T3}};
\node at (-3.5,-4.75) {\normalsize\textbf{T4}};
\node at (-2.5,-5.25) {\normalsize\textbf{T5}};
\node (v68) at (-7.5,-2.25) {};
\node (v69) at (5.25,-2.25) {};
\draw [fill = green, thick] (v68) edge (v69);
\node at (-1.25,-2.5) {Below: Overall airline schedule for multiple days.};
\node at (-1.25,-2) {Above: Detailed disruption management during tracking period for single flight day};
\end{tikzpicture}
\caption{Illustrative example detailing the practical context of the RCRP }
\label{fig:context}
\end{figure}

In Figure \ref{fig:context}, we provide a scheme that details the airline scheduling operations which facet the RCRP.  Five flight days are denoted (D1 - D5) and each of these flight days has a corresponding crew-schedule (both regular and reserves). Long before the actual flight day, an initial schedule is published (T1 - T5). In this paper, we focus on a single flight day with accompanying tracking period, as is denoted with the black ellipse. We consider that at the scheduling time point, the number of the available crew members are known. The figure's upper-part shows such a tracking period in detail. We indicate disruptions with the red arrows. When such a disruption occurs, we recover from it by solving the RCRP. This provides us with a recovered schedule for the flight day, for both regular crew and reserve crew. Every time a disruption occurs, the RCRP is solved which provides us with updated (reserve) crew schedules. Since the RCRP schedules reserve crew in a reliable way, as will be detailed in Section \ref{sec:reserves}, we end up with a (reserve) crew schedule that can effectively cope with disruptions at the day of execution.

\subsection{Problem description}

In the following, we formally introduce the RCRP. 

Let $\cal{F} = \{1, \ldots, F\}$ be the set of flight legs (in short, flights) in the schedule on the day under consideration. For each flight $f\in \cal{F}$, we consider a discrete set of flight copies $\Omega_f$, consisting of copies of flight $f$ with different departure times that are feasible with respect to the total airline schedule. We let $\Omega = \cup_{f \in \cal{F}}\Omega_f$, and we denote a single flight-copy as $\omega \in \Omega_f$. Whereas two flight copies $\omega, \omega' \in \Omega$ might be feasible on its own, selecting both might be infeasible due to, for instance, airport restrictions on flights departing and arriving at similar time slots or flights sharing the same aircraft. We let $\cal{I}:=\{(\omega, \omega') \mid \omega, \omega' \in \Omega \text{ and } \omega, \omega' \text{ are incompatible} \}$ be the set of all incompatible flight copies. 

Let $\cal{K} = \{1, \ldots, K\}$ be the set of all crew members on the day under consideration, of which $\cal{K}^R \subset K$ denotes the set of reserve-crew members and $\cal{K}^O \subset K$ denotes the set of regular crew-members.
For each crew member $k \in \cal{K}$, let $\cal{P}^k$ denote the set of possible crew pairings, and let $\cal{P} = \cup_{k\in\cal{K}}\cal{P}^k$. 

Five criteria should be met for a crew-pairing in order to comply with the \textit{duty legality rules}. First, the maximum time spent flying during a duty by a crew member is at most $\phi_1^\textsc{l}$  for a long duty (more than 4 flight legs) and $\phi_1^\textsc{s}$ for a short duty (less than 3 flight legs), respectively. Second, the maximum length (in hours) of a duty is at most $\phi_2^\textsc{L}$ and $\phi_2^\textsc{S}$ for a long and short duty, respectively. Third, the minimum time between two consecutive flights (sit time/turnaround time) should be at least $\phi_3$. Fourth, the minimum time between two duties (rest time) should be at least $\phi_4$. {}{In practice, rest time depends on the length of the duty. For the brevity of the model representation we use a single parameter to reflect this.} Lastly, the adjusted duties should fit in the overall flight schedule for an individual crew member. That is, crew members start their duties at the location where their previous duties end, and vice-versa.

Underlying the regular crew pairings that constitute a feasible flight plan, a reserve-crew schedule is present to mitigate risks from possible disruptions. Let $\mathcal{S} = \{1, \ldots, S\}$ be the set of reserve shifts. Each reserve shift $s \in \cal{S}$ has a fixed start time $B_s$ and a fixed end time $E_s$. We consider a set of reserve-crew schedules $\Theta_n$, where $n$ means the total number of reserve-crew members available. 
We define $\theta:=(h^\theta_1, h^\theta_2, \ldots, h^\theta_S)$ and $\theta \in \Theta_n$.
$\theta$ represents a reserve-crew schedule consisting of a sequence of numbers of crew members assigned to each shift with $\sum_{i=1}^S h^\theta_i = n$. That is, $h_s^\theta$ denotes the number of crew members assigned to shift $s \in \cal{S}$ in reserve-crew schedule $\theta \in \Theta_n$. Let $\Theta_n^k \subset \Theta_n$ be the set of reserve-crew schedules that contains crew member $k$. Note that shifts typically overlap and selecting which reserve-crew member to use is non-trivial. This will be explained in Section 3.2.

We consider a set $\cal{D} = \{1, \ldots, D\}$ of crew-induced disruptions, where each $d \in \cal{D}$ describes a flight leg no longer being covered. The goal is to recover (i.e., once again cover all flights) at minimum cost while retaining the flexibility in the resulting schedule to cope with potential additional disruptions in the future. In addition, one may choose to let crew deadhead (i.e., travel as a passenger) on a flight at the price of incurring deadhead costs $c^{d}_w, w \in \Omega$, or a flight might be canceled at the price of incurring  a flight copy dependent canceling cost $c^\text{c}_f$. 

\subsection{Mixed integer programming formulation} \label{sec:masterproblem}
Let $x_p^k$ be a binary decision variable that equals 1 if pairing $p$ is selected for crew member $k$, and equals 0 otherwise.  Let $a_{\omega p}^{k}$ be a binary parameter that is 1 if pairing $p$ of crew member $k$ includes flight copy $\omega$, and is 0 otherwise. Furthermore, we introduce artificial variables $y_{\omega}$  that denote the number of crew deadheading on flight copy $\omega \in \Omega$, and define binary variables $z_f$ being equal to 1 if flight $f$ is cancelled and zero otherwise. Let $u_\theta$ be a binary variable that equals 1 if reserve-crew schedule $\theta \in \Theta_n$ is selected, and is zero otherwise. The costs for selecting reserve-crew schedule $\theta$ equals $c_\theta$. Let $c_p^k$ be the costs of selecting pairing $p$ for crew-member $k$. We detail those costs in Section \ref{sec:coststructure}.

We condition on the number of available reserve-crew members, and we write RCRP($n$) to denote an RCRP instance with $n$ reserve-crew members. The RCRP($n$) is formulated as the following mixed-integer linear programming model. The list of notation is given in Table 1. 
\begin{align}
\min \quad  & \sum_{k \in \cal{K}} \sum_{p\in \cal{P}^k}c_p^k x_p^k + \sum_{\omega \in \Omega}{}{c^{d}_{\omega}}y_{\omega} + \sum_{f\in \cal{F}}c^c_f z_f + \sum_{\theta \in \Theta_n}c_\theta u_\theta,  \hspace{4.6cm}\label{lp:objective} 
\end{align} 
\vspace{-0.8cm}
\begin{align}
 \text{s.t. } \quad & \sum_{k\in \cal{K}} \sum_{p \in \cal{P}^k} \sum_{\omega \in \Omega_f}a^{k}_{\omega p} x_p^k - \sum_{\omega\in\Omega_f}y_{\omega} + z_f = 1 &\ & \quad \forall f\in \cal{F}, \label{lp:coveringcrew}\\
 &  y_{\omega} - \sum_{k \in \cal{K}} \sum_{p \in \cal{P}^k} a_{\omega p}^kx_p^k \leq 0   &\ & \quad \forall \omega \in \Omega \label{lp:deadheadconsistency1} \\
 & \sum_{w \in \Omega_f} y_{\omega} \leq M(1-z_f)   &\ & \quad \forall f \in \cal{F} \label{lp:deadheadconsistency3} \\
&  \sum_{k\in \cal{K}}\sum_{p\in \cal{P}^k}(a^{k}_{\omega p}x_p^k - y_{\omega}) +  \sum_{k\in \cal{K}}\sum_{p\in \cal{P}^k}(a^{k}_{\omega' p}x_p^k  - y_{\omega'}) \le 1&\ & \quad \forall (\omega, \omega') \in \cal{I}, \label{lp:illegalcopy}   \\
& \sum_{p \in \cal{P}^k} x_p^k  \leq  1 &\ & \quad \forall k\in \cal{K}^O, \label{lp:assignmentcrew} \\
& \sum_{p \in \cal{P}^k} x_p^k \le 1 &\ & \quad \forall k\in \cal{K}^R, \label{lp:assignmentreservecrew} \\
& \sum_{p\in \cal{P}^k} x_p^k - \sum_{\theta \in \Theta_n^k}{}{\lambda^{kp}_{\theta}} u_\theta \le 0  &\ & \quad \forall k \in \cal{K}^R, \label{lp:scheduleconsistency} \\
& \sum_{\theta \in \Theta_n} u_{\theta} \le 1  &\ & \quad \label{lp:scheduleassignment} \\
& x_{p}^k \in \{0,1\} &\ & \quad \forall k\in \cal{K},  \forall p\in \cal{P}^k, \label{lp:binarycrew} \\
& y_{\omega} \ge 0 &\ & \quad \forall  \omega \in \Omega, \label{lp:binarydeadhead1}  \\
& z_f \in \{0,1\} &\ & \quad \forall f\in \cal{F},
\label{lp:binarycancel}  \\
& u_\theta \in \{0,1\} &\ & \quad \forall \theta \in \Theta_n. \label{lp:binaryschedule} 
\end{align}

\begin{table} 
\caption{Summary of main notation}
\begin{tabular}{lp{4in}} 
\toprule 
	\textbf{Sets:} \\
	$\cal{F}$ & Set of flight legs \\
	$\cal{K}$ & Set of all crew members on the day under consideration \\
	$\cal{K}^R$  & Set of reserve crew members\\
	$\cal{K}^O$ & Set of regular crew members\\
	 $\cal{P}^k$ & Set of possible crew pairings including crew member $k$ \\
	 $S$ &  Set of reserve shifts\\
	  $\Omega$ & Set of flight copies \\
	  $\cal{I}$ & Set of all incompatible flight copies\\
	  	 $\Theta_n$ & Set of reserve-crew schedule with $n$ reserve crew members available \\
	\textbf{Parameters:} \\
	$a_{\omega p}^{k}$ & Binary parameter that is 1 if pairing $p$ of crew member $k$ includes flight copy $\omega$, and is 0 otherwise \\
$c_\theta$ &  The costs for selecting reserve-crew schedule $\theta$ \\
$c_p^k$ &  The cost of selecting pairing $p$ for crew-member $k$ \\
$c_{\omega}^d$ & The deadhead cost per crew member for flight leg $d$ and flight copy $\omega$\\
$c_{f}^c$ & The canceling cost for flight $f$\\
 $c_{t}$ & Fixed transportation costs between bases\\
 $c_\textsc{a}$ & Schedule altering costs\\
 $\lambda^{k p}_{\theta }$ &  Binary parameter that is 1 if a pairing $p$ can be matched by crew member $k$ in schedule $\theta$  \\
 $N^{\theta}$ & Maximum number of available reserve-crew members\\

	\textbf{Decision variables:} \\
	$x_{p}^k$ & Binary variable that equals 1 if pairing $p$ is selected for crew member $k$, and equals 0 otherwise\\
	$y_{\omega}$ & The number of crew deadheading on flight copy $\omega \in \Omega$ \\
    $z_f$ & Binary variables $z_f$ being equal to 1 if flight $f$ is cancelled and zero otherwise\\
    $u_\theta$ & Binary variable that equals 1 if reserve-crew schedule $\theta \in \Theta_n$ is selected, and is zero otherwise \\ \bottomrule
 	\end{tabular}
\end{table}

The Objective \eqref{lp:objective} consists of the costs for the selected crew pairings, the costs incurred for the number of deadheads, the costs incurred through canceling flights, and the costs incurred for the selected reserve schedule (and thus, for the selection which reserves to use). Constraints (\ref{lp:coveringcrew}) models (in combination with Constraints (3) and (4)) that if we cancel a flight, we do not allow deadheading (i.e., $y_w$) on any flight copy and we do not allow to assign a crew member to that flight (i.e., $a_{\omega p}^k x_p^k$). If we do not cancel the flight, we should always assign a crew member to that flight and we let crew members deadhead on that flight. {Note that a crew member can be interpreted both as an individual team member and as an entire team, all of whom are subject to the same work regulations. This is inspired upon the collaboration with our industry partner.}
Constraints (\ref{lp:deadheadconsistency1}) and (\ref{lp:deadheadconsistency3}) ensure that no crew is deadheading on a flight-copy which is not operated by at least one active crew member, so that consistency of departure times between crew flying and crew deadheading is guaranteed. 
Constraints \eqref{lp:illegalcopy} ensure that no illegal pairs of flight copies are used, and  Constraints (\ref{lp:assignmentcrew}) ensure that not-disrupted regular crew  is assigned to exactly one pairing. Reserve crew members are restricted to at most one pairing through Constraints (\ref{lp:assignmentreservecrew}). Crew members without a pairing selected, start and end at the location where they originated from at the beginning of the day. {}{ Constraints (\ref{lp:scheduleconsistency}) ensure  consistency of the reserve schedule accompanying the selected reserve crew members where $\lambda^{kp}_{\theta}$ indicates whether a pairing $p$
can be matched by crew member $k$ in schedule $\theta$.} Finally, Constraints  (\ref{lp:scheduleassignment}) ensure that at most one reserve schedule is selected. Note that the formulation is as such that feasibility is always guaranteed (i.e., by canceling all the flights).

The set of pairings  $p \in \cal{P}^k$ as encountered in the RCRP($n$) is cumbersome to enumerate up-front, and therefore column generation is an appropriate method to solve the LP relaxation of RCRP($n$). We will refer to the linear relaxation of RCRP($n$) as L-RCRP($n$), where we replace the integrality restrictions of the binary decision variables with $x_p^k \geq 0$ for all $p \in \cal{P}^k$ and $k \in \cal{K}, z_f \geq 0$ for all $f \in \cal{F}$, and $u_\theta \geq 0$ for all $\theta \in \Theta$. To solve L-RCRP($n$) with column generation, we will work with a restricted set of pairings $\overline{\cal{P}}^k \subset \cal{P}^k$, and we will call the L-RCRP($n$) subject to these restricted set of pairings, the restricted linear relaxation of RCRP($n$) (RL-RCRP($n$)). Then, column generation can be described as the procedure which iteratively: 1) solves RL-RCRP($n$), and 2) generates for each $k \in \cal{K}$ pairings $p \in \cal{P}^k \backslash \overline{\cal{P}}^k$ with negative reduced cost by solving a so-called pricing problem (defined below), and includes those pairings in $\overline{\cal{P}}^k$. Optimality is ensured if no pairings of negative reduced cost are found after solving the pricing problem, in other words, the corresponding dual solution of RL-RCRP($n$) is feasible for L-RCRP($n$), and thereby optimal for L-RCRP($n$) as well. 

In order to formulate the pricing problem, let us define the appropriate dual variables. Note that the variables may depend on an arbitrary crew member $k$ due to decomposing the problem in the crew-member dimension. Let  $\alpha_f$ and $\beta_{\omega}$ be the dual variables corresponding to constraints \eqref{lp:coveringcrew} and \eqref{lp:deadheadconsistency1}, respectively. Dual variables for the infeasible flight-copy constraints \eqref{lp:illegalcopy} are given by $\gamma_{\omega,\omega'}$. Let $\delta^k$ be the dual variables corresponding to Constraints \eqref{lp:assignmentcrew} and \eqref{lp:assignmentreservecrew}, and $\epsilon^k$ be the dual variables corresponding to Constraints \eqref{lp:scheduleconsistency}. Then the pricing problems for a regular crew member $k \in \cal {K}^O$ asks for solving
\begin{equation}
    \min_{p\in\cal{P}^k} \hat{c}_p^k := c_p^k - \sum_{f\in \cal{F}}\sum_{\omega \in \Omega_f}\alpha_f a_{\omega p}^k - \sum_{\omega \in \Omega}\beta_{\omega} a^{k}_{\omega p} - \sum_{(\omega,\omega') \in \cal{I}}\gamma_{\omega,\omega'}(a^{k}_{\omega p} + a^{k}_{\omega' p}) - \delta^k,\label{eq:p1}
\end{equation}
and the pricing problem for a reserve-crew member $k' \in \cal{K}^R$ asks for solving
\begin{equation}
       \min_{p\in\cal{P}^{k'}} \hat{c}_p^{k'} := c_p^{k'} - \sum_{f\in \cal{F}}\sum_{\omega \in \Omega_f}\alpha_f a_{\omega p}^{k'} - \sum_{\omega \in \Omega} \beta_{\omega} a^{k'}_{\omega p} - \sum_{(\omega,\omega') \in \cal{I}}\gamma_{\omega,\omega'}(a^{k'}_{\omega p} + a^{k'}_{\omega' p}) - \delta^{k'} - \epsilon^{k'}. \label{eq:p2}
\end{equation}
The pricing problems depicted above are variants of the Resource Constrained Shortest Path Problem \citep[see, e.g.,][]{irnich2005}. We will further elaborate on the structure of the pricing problems in Section  \ref{sec:pricingproblem}, where we provide an efficient algorithm to solve the pricing problems.

\section{Cost structures} \label{sec:coststructure}
In this section, we introduce the two parts of the RCRP($n$) that are still undefined. Namely, we detail the crew paring costs $c_p^k$ (Section \ref{sec:pairingcost}) and the reliable reserve crew schedule costs $c_\theta$ (Section \ref{sec:reserves}). The RCRP selects a complete reserve crew schedule covering all shifts (and thus periods) for the flight execution day. Thus, the cost $c_{\theta}$ for selecting a reserve crew schedule (i.e., the expected shortfall) completely covers all periods and shifts on the day of execution.

\subsection{Pairing costs for regular crew}\label{sec:pairingcost}
The costs of a crew pairing $c^k_p$ (i.e., an assignment of a crew member to a duty) consist of several parts, that are detailed next. Following \cite{barnhart2003b}, we will consider the pay-and-credit for a pairing p, referred to as $PC(p)$. For any crew-pairing $p$, let $FT(p)$ denote the cost associated with the flight time of the pairing and let $D(p)$ be the cost related to the total duration of the pairing. Let $PC^{\min}$ be the minimum guaranteed pay-and-credit for any pairing. The pay-and-credit $PC(p)$ is then defined as
\begin{equation}
    PC(p) = \max\{\textit{FT(p)}, D(p), PC^{\min} \}.
\end{equation}

In addition to the pay-and-credit costs for a pairing $p$, transportation costs arise whenever the origin and destination of a crew member's pairing are at different locations. These costs are denoted by $TC(p)$ for each pairing $p \in \cal{P}$, where {}{$TC(p) := I_{\{ k(p) \in \cal{K}^O\}}c_t$. Here, $k(p)$ denotes the crew member corresponding to pairing $p$, $\cal{K}^O$ is the set of the regular crew members, and $c_t$ is the fixed transportation cost between bases for one crew member. $I$ is the indicator function showing whether a crew member is in a pair $p$.  } These transportation costs are defined as a symmetric cost function taking as arguments the origin and destination location.

Evaluating a pairing solely by its planned costs without disruptions ignores the obvious effect of the pairing structure on the robustness of the overall crew schedule. We, therefore, consider explicit penalties on the following four unfavorable characteristics of a pairing $p$.
\begin{enumerate}
    \item  The pairing's flight time is close to the maximum allowed flight time $\phi_1^L$ or $\phi_1^S$.
    \item  The pairing's duration is close to the maximum allowed pairing duration $\phi_2^L$ or $\phi_2^S$.
    \item  The pairing's sit time between consecutive flight legs is close to the minimum required rest time $\phi_3$.
    \item  The rest time between consecutive flight legs is close to the minimum required time between consecutive flight legs $\phi_4$.
\end{enumerate}
The total penalty costs $PEN(p)$ of pairing $p$ is then given by
\begin{equation}
    PEN(p) = \sum_{i=1}^4\left(g_i - \Delta_i\abs{\chi_i(p) - \phi_i}\right)^+,
\end{equation}
where $g_i$ equals the maximum penalty for characteristic $i$,  $\Delta_i$  is a parameter that determines the slope with which the penalty decreases, $\chi_i(p)$ is the observed value for characteristic $i$ in pairing $p$, and $\phi_i$ equals the maximum or minimum value for the given characteristic. 

By equation (18), the penalty for characteristic $i$ is the maximum penalty $g_i$ iff $ \chi_i(p) = \phi_i$. This penalty takes into account uncertainty and risk of being beyond the required threshold for the corresponding characteristic before the reliable modelling itself.

Re-timing costs $c_f^\textsc{r}$ are included when a flight included  pairing is re-timed. Let $RET(p)$ be the sum of those costs for pairing $p$, i.e,. $RET(p):= \sum_{f\in F}\sum_{\omega \in \Omega_f}a^{k}_{\omega p} c_f^\textsc{r}$. Note that we deduct these costs from the deadheading costs $c^d_{\omega}$ for non-original flight copies $\omega \in \Omega_f$, otherwise the pairing associated with deadheading would also be penalized for being re-timed. {}{By deducting the cost from the deadhead cost, we ensure that if a crew member deadheads on a re-timed pairing, it will only incur deadheading cost and not re-timing cost.}

Finally, altering costs are incurred if a crew member gets assigned a pairing different than the already published schedule prescribes. Let $ALT(p) := I_{\{ k(p) \in \cal{K}^O\}}c_a$ be those costs, where $c_a$ means the scheduling altering cost for one crew member. Note that the altering costs will equal zero for reserve-crew members, and equal $c_a$ for regular crew members.

Then, the total pairing costs for a pairing $p$ are defined as
\begin{equation}
    c^k_p = \max\{PC(p) + TC(p) + PEN(p) - c_\textsc{o}^k, 0\} + RET(p) + ALT(p),
\end{equation}
where $c_\textsc{o}^k$ are the costs of the original pairing of member $k$, i.e., the costs of its corresponding pairing in the already published schedules. {}{Since we intend to make the optimal adjustment based on the already published schedules, we only counter the potential cost increase caused by the adjustment.}

\subsection{The reliable reserve-crew schedule costs} \label{sec:reserves}
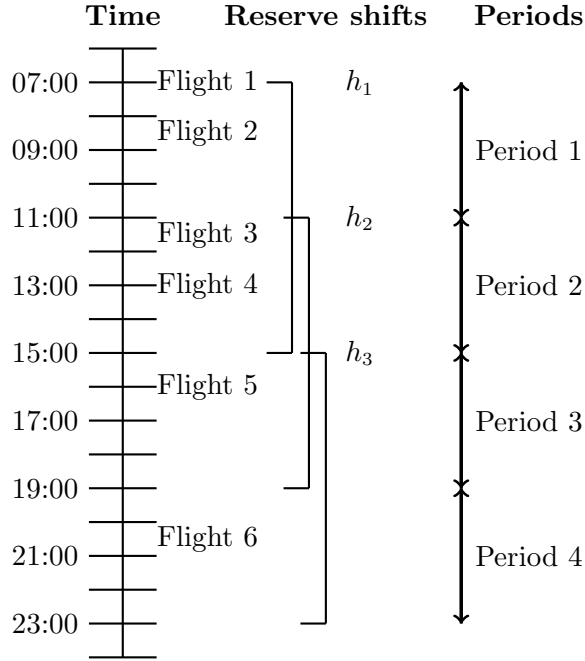
\begin{figure}[h!]
\centering
\begin{tikzpicture}[scale = 0.45]
    \newcommand{\width}{1}
    \newcommand{\height}{1}

    \coordinate (A1)         at     (0,0) {};
    \coordinate (A2)         at     ($(A1)+(-1*\width, 0)$) {};
    \coordinate (A3)         at     ($(A2)+(-1*\width, 0)$) {};
    \coordinate (B1)         at     ($(A1)+(0, 1*\height)$) {};
    \coordinate (B2)         at     ($(B1)+(-1*\width, 0)$) {};
    \coordinate (B3)         at     ($(B2)+(-1*\width, 0)$) {};
    \coordinate (C1)         at     ($(B1)+(0, 1*\height)$) {};
    \coordinate (C2)         at     ($(C1)+(-1*\width, 0)$) {};
    \coordinate (C3)         at     ($(C2)+(-1*\width, 0)$) {};
    \coordinate (D1)         at     ($(C1)+(0, 1*\height)$) {};
    \coordinate (D2)         at     ($(D1)+(-1*\width, 0)$) {};
    \coordinate (D3)         at     ($(D2)+(-1*\width, 0)$) {};
    \coordinate (E1)         at     ($(D1)+(0, 1*\height)$) {};
    \coordinate (E2)         at     ($(E1)+(-1*\width, 0)$) {};
    \coordinate (E3)         at     ($(E2)+(-1*\width, 0)$) {};
    \coordinate (F1)         at     ($(E1)+(0, 1*\height)$) {};
    \coordinate (F2)         at     ($(F1)+(-1*\width, 0)$) {};
    \coordinate (F3)         at     ($(F2)+(-1*\width, 0)$) {};
    \coordinate (G1)         at     ($(F1)+(0, 1*\height)$) {};
    \coordinate (G2)         at     ($(G1)+(-1*\width, 0)$) {};
    \coordinate (G3)         at     ($(G2)+(-1*\width, 0)$) {};
    \coordinate (H1)         at     ($(G1)+(0, 1*\height)$) {};
    \coordinate (H2)         at     ($(H1)+(-1*\width, 0)$) {};
    \coordinate (H3)         at     ($(H2)+(-1*\width, 0)$) {};
    \coordinate (I1)         at     ($(H1)+(0, 1*\height)$) {};
    \coordinate (I2)         at     ($(I1)+(-1*\width, 0)$) {};
    \coordinate (I3)         at     ($(I2)+(-1*\width, 0)$) {};
    \coordinate (J1)         at     ($(I1)+(0, 1*\height)$) {};
    \coordinate (J2)         at     ($(J1)+(-1*\width, 0)$) {};
    \coordinate (J3)         at     ($(J2)+(-1*\width, 0)$) {};
    \coordinate (K1)         at     ($(J1)+(0, 1*\height)$) {};
    \coordinate (K2)         at     ($(K1)+(-1*\width, 0)$) {};
    \coordinate (K3)         at     ($(K2)+(-1*\width, 0)$) {};
    \coordinate (L1)         at     ($(K1)+(0, 1*\height)$) {};
    \coordinate (L2)         at     ($(L1)+(-1*\width, 0)$) {};
    \coordinate (L3)         at     ($(L2)+(-1*\width, 0)$) {};
    \coordinate (M1)         at     ($(L1)+(0, 1*\height)$) {};
    \coordinate (M2)         at     ($(M1)+(-1*\width, 0)$) {};
    \coordinate (M3)         at     ($(M2)+(-1*\width, 0)$) {};
    \coordinate (N1)         at     ($(M1)+(0, 1*\height)$) {};
    \coordinate (N2)         at     ($(N1)+(-1*\width, 0)$) {};
    \coordinate (N3)         at     ($(N2)+(-1*\width, 0)$) {};
    \coordinate (O1)         at    ($(N1)+(0, 1*\height)$) {};
    \coordinate (O2)         at     ($(O1)+(-1*\width, 0)$) {};
    \coordinate (O3)         at     ($(O2)+(-1*\width, 0)$) {};
    \coordinate (P1)         at     ($(O1)+(0, 1*\height)$) {};
    \coordinate (P2)         at     ($(P1)+(-1*\width, 0)$) {};
    \coordinate (P3)         at     ($(P2)+(-1*\width, 0)$) {};
    \coordinate (Q1)         at     ($(P1)+(0, 1*\height)$) {};
    \coordinate (Q2)         at     ($(Q1)+(-1*\width, 0)$) {};
    \coordinate (Q3)         at     ($(Q2)+(-1*\width, 0)$) {};
    \coordinate (R1)         at     ($(Q1)+(0, 1*\height)$) {};
    \coordinate (R2)         at     ($(R1)+(-1*\width, 0)$) {};
    \coordinate (R3)         at     ($(R2)+(-1*\width, 0)$) {};
    \coordinate (S1)         at     ($(R1)+(0, 1*\height)$) {};
    \coordinate (S2)         at     ($(S1)+(-1*\width, 0)$) {};
    \coordinate (S3)         at     ($(S2)+(-1*\width, 0)$) {};
    
    \draw[thick] (A2) -- (S2);
    
    \draw[thick] (A1) -- (A3) node[anchor=east] {};
    \draw[thick] (B1) -- (B3) node[anchor=east] {23:00};
    \draw[thick] (C1) -- (C3) node[anchor=east] {};
    \draw[thick] (D1) -- (D3) node[anchor=east] {21:00};
    \draw[thick] (E1) -- (E3) node[anchor=east] {};
    \draw[thick] (F1) -- (F3) node[anchor=east] {19:00};
    \draw[thick] (G1) -- (G3) node[anchor=east] {};
    \draw[thick] (H1) -- (H3) node[anchor=east] {17:00};
    \draw[thick] (I1) -- (I3) node[anchor=east] {};
    \draw[thick] (J1) -- (J3) node[anchor=east] {15:00};
    \draw[thick] (K1) -- (K3) node[anchor=east] {};
    \draw[thick] (L1) -- (L3) node[anchor=east] {13:00};
    \draw[thick] (M1) -- (M3) node[anchor=east] {};
    \draw[thick] (N1) -- (N3) node[anchor=east] {11:00};
    \draw[thick] (O1) -- (O3) node[anchor=east] {};
    \draw[thick] (P1) -- (P3) node[anchor=east] {09:00};
    \draw[thick] (Q1) -- (Q3) node[anchor=east] {};
    \draw[thick] (R1) -- (R3) node[anchor=east] {07:00};
    \draw[thick] (S1) -- (S3) node[anchor=east] {};

    \node[] at ($(R1)+(1.5*\width, 0)$)  {Flight 1};
    \node[] at ($(Q1)+(1.5*\width, -0.5*\height)$)  {Flight 2};
    \node[] at ($(N1)+(1.5*\width, -0.5*\height)$)  {Flight 3};
    \node[] at ($(L1)+(1.5*\width, 0)$)  {Flight 4};
    \node[] at ($(I1)+(1.5*\width, 0)$)  {Flight 5};
    \node[] at ($(E1)+(1.5*\width, -0.5*\height)$)  {Flight 6};
    
    \coordinate (Line1)         at     ($(R1)+(4*\width, 0)$) {};
    \coordinate (Line1b)        at     ($(Line1) + (-0.75*\width, 0)$) {};
    \coordinate (Line2)         at     ($(J1)+(4*\width, 0)$) {};
    \coordinate (Line2b)        at     ($(Line2) + (-0.75*\width, 0)$) {};
    \coordinate (Line3)         at     ($(N1)+(4.5*\width, 0)$) {};
    \coordinate (Line3b)        at     ($(Line3) + (-0.75*\width, 0)$) {};
    \coordinate (Line4)         at     ($(F1)+(4.5*\width, 0)$) {};
    \coordinate (Line4b)        at     ($(Line4) + (-0.75*\width, 0)$) {};
    \coordinate (Line5)         at     ($(J1)+(5*\width, 0)$) {};
    \coordinate (Line5b)        at     ($(Line5) + (-0.75*\width, 0)$) {};
    \coordinate (Line6)         at     ($(B1)+(5*\width, 0)$) {};
    \coordinate (Line6b)        at     ($(Line6) + (-0.75*\width, 0)$) {};
    
    \draw[thick] (Line1) -- (Line2);
    \draw[thick] (Line3) -- (Line4);
    \draw[thick] (Line5) -- (Line6);
    
     \draw[thick] (Line1) -- (Line1b);
     \draw[thick] (Line2) -- (Line2b);
     \draw[thick] (Line3) -- (Line3b);
     \draw[thick] (Line4) -- (Line4b);
     \draw[thick] (Line5) -- (Line5b);
     \draw[thick] (Line6) -- (Line6b);
     
      \coordinate (Point1)         at     ($(Line1)+(5*\width, 0)$) {};
      \coordinate (Point2)         at     ($(Line3)+(4.5*\width, 0)$) {};
      \coordinate (Point3)         at     ($(Line5)+(4*\width, 0)$) {};
      \coordinate (Point4)         at     ($(Line4)+(4.5*\width, 0)$) {};
      \coordinate (Point5)         at     ($(Line6)+(4*\width, 0)$) {};
      
      \node[] at ($(Point1)+(-3*\width, 0)$)  {$h_1$};
      \node[] at ($(Point2)+(-3*\width, 0)$)  {$h_2$};
      \node[] at ($(Point3)+(-3*\width, 0)$)  {$h_3$};
      
       \draw[very thick, ->] (Point1) -- (Point2);
       \draw[very thick, ->] (Point2) -- (Point1);
       \draw[very thick, ->] (Point2) -- (Point3);
       \draw[very thick, ->] (Point3) -- (Point2);
       \draw[very thick, ->] (Point3) -- (Point4);
       \draw[very thick, ->] (Point4) -- (Point3);
       \draw[very thick, ->] (Point4) -- (Point5);
       \draw[very thick, ->] (Point5) -- (Point4);
       
       \node[] at ($(Point1)+(2*\width, -2*\height)$)  {Period 1};
       \node[] at ($(Point2)+(2*\width, -2*\height)$)  {Period 2};
       \node[] at ($(Point3)+(2*\width, -2*\height)$)  {Period 3};
       \node[] at ($(Point4)+(2*\width, -2*\height)$)  {Period 4};
       
       \node[] at ($(S2)+(0, 1*\height)$)  {\textbf{Time}};
       \node[] at ($(Line3)+(0.5, 6*\height)$)  {\textbf{Reserve shifts}};
       \node[] at ($(Point1)+(2, 2*\height)$)  {\textbf{Periods}};
\end{tikzpicture}
\caption{Illustrative example of a reserve-crew schedule}
\label{fig:timeline}
\end{figure}

Recall that a reserve-crew schedule $\theta \in \Theta_n$ is defined as $\theta = (h_1^\theta, \ldots, h_S^\theta)$ as an assignment of $h_s^\theta$ crew-members to shifts $s \in \cal{S}$. In this section, we develop a robust cost-measure ($c_\theta$) for using reserve-crew schedule $\theta$. Recall the RCRP model focuses on the so-called flight execution day. In the following, we will develop a Markov chain model that models operations on that flight execution day in detail and use that to derive $c_\theta$. We start with an illustrative example where we introduce the additional required terminology.
\begin{example}
Figure \ref{fig:timeline}  depicts the time horizon for an arbitrary day, starting at 06:00 h and ending at 23:59 h. Six flights and three reserve shifts are scheduled. We assume that all reserve shifts are of length $\tau$. In Figure \ref{fig:timeline}, three reserve shifts are depicted with $h_1, h_2$, and $h_3$ crew members assigned, respectively. Based on the schedule shifts, we define the set $\cal{T}$ of periods within a day. Note that, shifts and periods are two different concepts. One shift may consist of multiple periods and may start even if the previous shift has not finished yet. 
In Figure \ref{fig:timeline}, four periods are defined, where, for instance, shift 2 consists of two periods starting at the beginning of period 2 and ending at the end of period 3. In order to capture all flights, it should be ensured that no flight occurs before the start of the first period or the end of the last period.
\end{example}

We model the number of reserve-crew members available at period $t \in \cal{T}$ as a Markov-chain with state space $\cal{E} = \{{}{0}, \ldots, N^\theta \}$, where $N^\theta$ is the maximum number of available reserve-crew members in reserve-crew schedule $\theta \in \Theta_n$. {}{The order of events in each period $t$ is as follows. First, we consider the reserve crew members that finished their shift at the start of period $t$ (or equivalently at the end of period $t-1$, but for notational convenience we consistently write that this happens at the start of period $t$). Second, we consider the reserve crew members that start their shift at the start of period $t$. Third, we incur the stochastic demand for reserve crew members during period $t$.}

Let $Q^\theta_t$ describe the distribution of remaining reserve crew {}{at the end of period $t$}, i.e.,  $Q_t^\theta = (\pi^\theta_{0t},\pi^\theta_{1t},\dots,\pi^\theta_{N^\theta t})$, where $\pi^\theta_{et}$ denotes the probability of having $e$ reserve crew members {}{ at the end of period t, or equivalently, at the start of period $t+1$ in reserve-crew schedule $\theta$ as no further events happen in period $t$ after observing $Q^\theta_t$.}

\subsubsection{An optimal recovery policy.}
To correctly model $Q_t^\theta$, a recovery policy describing how and when to use reserve crew will be defined in the following. In order to do so, we make the following three assumptions.

\begin{assumption}
Demand for reserve-crew members in period $t$ is described by a known discrete probability distribution $F_t$, with density $f_t$ and domain $[0, b_t]$. {}{Here, $F_t$ can for example be estimated based on historical data.}
\end{assumption} 

\begin{assumption}
Reserve crew demand is independently distributed among periods.  
\end{assumption}

\begin{assumption}
Reserve crew are always available to take over a flight, as long as it departs within their associated shift. 
\end{assumption}
{}{In the following, we consider three transition matrices that model transitions due to reserve crew members finishing their shift (step 1), starting their shift (step 2), and the stochastic demand for reserve crew members (step 3). These three transition matrices as a function of the period $t$ determines how $Q_{t-1}$ transitions to $Q_{t}$.}

\textbf{Step 1: Reserve crew ending their shift:} First, we will define transition matrix $P_t^{C\theta}$ that describes the transitions due to reserve crew ending their shift at the start of period $t$. We impose the Earliest-Finisher-First (EFF) recovery policy. This implies that if there is demand for a reserve-crew members in a period, we always assign the reserve-crew members that finish their shift earliest. The EEF policy allows us to calculate based on the observed number of remaining reserve-crew members (i.e., the ones that are not used to fulfill demand) how many will end their shift in the transition associated with step 1. Note that this policy corresponds to a First-In-First-Out policy if the shift lengths $\tau$ are equal for all reserve crew. We will prove that this recovery policy is optimal under \mbox{Assumptions 1 - 3}.

\begin{proposition}
When maximizing the number of reserve crew available at any point in time, an EFF recovery policy is optimal. 
\end{proposition}
\proof{Proof.}
Let an arbitrarily reserve-crew schedule $\theta$ be given. Let $(s,t)$ be an arbitrary state, where $s$ is the number of reserve-crew members available at period $t$. Let $\pi_{et'}^{\theta}$ be the probability of having $e$ reserve-crew members after period $t' \geq t$ in reserve-crew schedule $\theta$ according to the EFF policy. We show that for any other policy $\tilde{\pi}$, it will hold that $\pi_{et'}^{\theta} \geq \tilde{\pi}_{et'}^{\theta}$ for all $t \geq t'$. 

Suppose we follow policy $\pi$. That means that at state $(s, t)$, upon scheduling a reserve-crew member, we select the crew member $x$ whose shift ends earliest. By construction, this leaves the maximum number of crew members in all following periods. Thus, $\pi_{et'}^{\theta} \geq \tilde{\pi}_{et'}^{\theta}$ for all $t \geq t'$ for any policy $\tilde{\pi}$ unequal to the EFF. $\square$ \endproof

{}{Let the start of period $t$ coincide with an end point $E_s$ of shift $s$, such that the reserve-crew members scheduled in shift $s$ are off duty at the start of period $t$. Note if a period does not coincide with the end of a shift the transition associated with step 1 is redundant and $P_t^{C\theta}$ is the identity matrix. Define $\bar{e}^\theta_{t-1} := \max \{e \mid \pi^\theta_{e,t-1} > 0\}$ as the maximum number of reserve crew members that can be available at the end of period $t-1$.} Assuming the EFF recovery policy, the transition matrix $P_t^{C\theta}$ is given by:
\begin{equation}
\label{eq22}
   (P_t^{C\theta})_{ij} = 
    \begin{cases}
    1 & \text{if } i > j \text{ and } j = \bar{e}_{t-1}^\theta - h_s^\theta \\ 
    1 & \text{if } i = j \text{ and } j \leq \bar{e}_{t-1}^\theta - h^\theta_s, \\
    0 & \text{otherwise}.
    \end{cases}
\end{equation}
The rationale behind each of the cases in Equation \ref{eq22} is as follows. The first case considers the case in which some reserve crew members of shift $s$ will end their duty \textit{and} these reserve crew members are not all being used to fulfill demand for reserve-crew members. This happens if the observed number of reserve-crew members is larger than the maximum number of reserve-crew members at $t-1$ ($\bar{e}_{t-1}^\theta)$ minus the number of reserve crew member ending their shift ($h_s^\theta$) at the start of period $t$. Hence, we transition with probability 1 from $i$ to $\bar{e}^\theta_{t-1} - h^{\theta}_s$ reserve crew members if $i > \bar{e}^\theta_{t-1} - h^{\theta}_s$. The second case covers all the situations in which we observe less than $\bar{e}^\theta_{t-1} - h^{\theta}_s$ reserve-crew members, implying that all the reserve crew members that end their shift at the start of period $t$ are already used to fulfill the demand for reserve-crew members. The third case means that the rest of state transitions cannot occur.

\textbf{Step 2: New reserve crew members starting their shift:} Let the transition matrix $P_t^{A\theta}$ describe reserve crew members starting their shift at the start of period $t$. Recall from Assumption $3$ that reserve crew are always available. Given $h_s'^{\theta}$ reserve crew starting their shift $s'$ at period $t$, the elements of $P_t^{A\theta}$ are given as
\begin{equation}
    (P_t^{A\theta})_{ij} = 
    \begin{cases}
    1 & \text{if } if i + h_{s'}^\theta \leq N^\theta \text{and} j = i + h^\theta_{s'} \\
    1 & {}{\text{if } i + h_{s'}^\theta > N^\theta \text{ and } j = i} \\
    0 & \text{otherwise}.
    \end{cases}
\end{equation}
Here, the first case is trivial as it increases the number of reserve crew members with the number of reserve crew members starting their shift $h_s^\theta$ if the observed number of crew members is a feasible observation. The second case is a technicality to ensure that the transition matrices are of the same dimension.

\textbf{Step 3: Reserve crew assigned to duties:} With the independent distribution of reserve-crew demand, the transition matrices due to reserves being assigned to a duty are described by $P_t^{B\theta}$, where the elements of $P_t^{B\theta}$ are given as
\begin{equation}
    (P_t^{B\theta})_{ij} = 
    \begin{cases}
    f_{t}(x) & \text{if } j = i-x > 0, \\ 
    1 - F_{t}(i) & \text{if } j = 0, \\
    0 & \text{otherwise},
    \end{cases}
\end{equation}
where $x$ means the demand of reserve-crew members.

\begin{figure}[!h]
\centering
\begin{tikzpicture}[auto, -triangle 60, semithick]
  \tikzstyle{every state}=[auto, draw=black, minimum width=0.8cm,]
   \tiny 

 \node[state]      (2_1)        
 {$2$};
  \node[state]      (2_2)[right = 2cm of 2_1]                   {$2$};
  \node[state]      (2_3)[right = 2cm of 2_2]                   {$2$};   
  \node[state]      (2_4)[right = 2cm of 2_3]                   {$2$};
  \node[state]      (1_1)[below = 1cm of 2_1]                  {$1$};
  \node[state]      (1_2)[right = 2cm of 1_1]                   {$1$};
  \node[state]      (1_3)[right = 2cm of 1_2]                   {$1$};   
  \node[state]      (1_4)[right = 2cm of 1_3]                   {$1$};
  \node[state]      (0_1)[below = 1cm of 1_1]                  {$0$};
  \node[state]      (0_2)[right = 2cm of 0_1]                   {$0$};
  \node[state]      (0_3)[right = 2cm of 0_2]                   {$0$};   
  \node[state]      (0_4)[right = 2cm of 0_3]                   {$0$};
  \node[state, scale = 0.01] (T_1)[above = 1cm of 2_3]                      {};

  \path
        (2_1) edge       node[sloped,pos=0.5] {$1$}  (1_2)
        (1_1) edge       node[sloped,pos=0.5] {$1$} (1_2)
        (0_1) edge       node[sloped,pos=0.5] {$1$}  (0_2)
        (2_2) edge       node[sloped,pos=0.7] {$1$}  (2_3)
        (1_2) edge       node[sloped,pos=0.7] {$1$} (2_3)
        (0_2) edge       node[sloped,pos=0.7] {$1$}  (1_3)
       (2_3) edge       node[sloped,pos=0.2] {$f_t(0)$} (2_4)
       (2_3) edge       node[sloped,pos=0.35] {$f_t(1)$}  (1_4)
       (2_3) edge       node[sloped,pos=0.65] {$1 - F_t(2)$}  (0_4)
      (1_3) edge       node[sloped,pos=0.2] {$f_t(0)$} (1_4)
        (1_3) edge       node[sloped,pos=0.05] {$1 - F_t(1)$} (0_4)
       (0_3) edge       node[sloped,pos=0.35] {$1$}  (0_4)
        
        ;
\draw[decoration={brace,mirror,raise=20pt},-,decorate, ]
   (0_1) -- node[below = 25pt] {\small $P^{C\theta}_t$} (0_2);
\draw[decoration={brace,mirror,raise=20pt},-,decorate]
   (0_2) -- node[below = 25pt] {\small $P^{A\theta}_t$} (0_3);
\draw[decoration={brace,mirror,raise=20pt},-,decorate]
   (0_3) -- node[below = 25pt] {\small $P^{B\theta}_t$} (0_4);

\end{tikzpicture}
\caption{Example of the state space transitioning with transition probabilities on the edges}
\label{fig:statespace}
\end{figure}
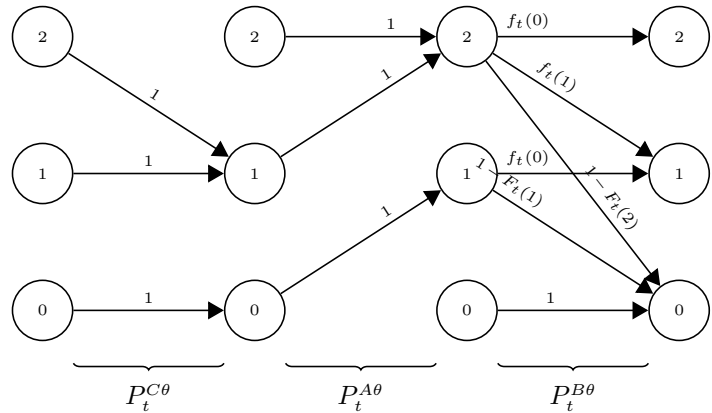

Then, the total transition matrix from period {}{$t-1$ to period $t$ is given by} $P_t^{C\theta}P_t^{A\theta}P_t^{B\theta}$, such that
\begin{equation}
    {}{Q^\theta_{t} = Q^\theta_{t-1}P_t^{C\theta}P_t^{A\theta}P_t^{B\theta}}.
\end{equation}
An illustrative example of such a transition is provided in Example \ref{ex:states}.
\begin{example}\label{ex:states}
In Figure \ref{fig:statespace}, we illustrate the transition matrices $P^{A\theta}, P^{B\theta}$, and $P^{C\theta}$ for a transition resulting from $1$ reserve crew ending his shift and $1$ reserve crew starting his shift with $ N^\theta= \bar{e}_{t-1}^\theta = 2$. From left to right, we indicate the transition matrices $P_t^{C\theta}$, $P_t^{A\theta}$, and $P_t^{B\theta}$. Since $\bar{e}_{t-1}^\theta = 2$ and 1 crew member ends his shift in this particular example, we know that if we observe 2 crew-members at the start of this transition (the upper left node), we will move to a state of 1 crew-member resulting from the crew-member ending his shift. {}{However, if we observe 0 or 1 crew member, we are sure that the crew-member that ends his shift is used to fulfill demand for reserve-crew members in some period $t' < t$. Hence, if we observe 0 or 1 crew member we are sure these crew members do not end their shift, and thus the number of reserve-crew remains equal. This explains the horizontal arcs from the nodes on the left. The transition matrix  $P_t^{A\theta}$ in the middle shows that 1 new crew-member starts his duty at this point. Note that arrow from 2 crew-members to 2 crew-members in the transition of $P_{t}^{A\theta}$ is the technicality mentioned earlier. The right transitions $P_t^{B\theta}$ show the probability of using crew members during period $t$ caused the demand for reserve-crew members induces by disruptions to the actual flight schedule.}
\end{example}

\endproof

\subsubsection{The risk of a reserve crew schedule.}
Given the transition matrices $P^{A\theta}_t, P^{B\theta}_t$, and $P^{C\theta}_t$ as described formerly, we are able to evaluate the risks associated with a given reserve crew schedule $\theta = (h^\theta_1, \dots, h^\theta_{S})$. We express these risks in terms of the expected shortfall (also called the conditional value at risk), which is defined as the expected shortage of reserve-crew members given that there is a shortage. Let $\xi_t(\theta)$ be the expected shortfall of reserve personnel in period $t$ for schedule $\theta$, i.e.,
\begin{equation}
\xi_t(\theta) = \sum_{i=1}^{b_t}i\cdot Q^\theta_{t-1} P_t^{C\theta} P_t^{A\theta} \cdot \mathbf{g^i_t}
\end{equation}
{}{where $b_t$ denotes the upper bound on reserve demand in period $t$ and $\mathbf{g^i_t}$ is a row vector of size $N^\theta + 1$ obtained from the probability density function $f_t$ for demand for reserve crew in period $t$ by setting $g^i_t[a] = f_t(a+i)$ if $a+i \leq b_t$, and 0 otherwise. The equation conditions over having $i$ crew members short (the summation over $i$). For each $i$, we multiply the state probability vector with $g^i_t$ to evaluate each state with a demand realization resulting in $i$ crew members being short. In this way, we evaluate each potential realization of demand for reserve crew personnel in period $t$, which can result from multiple realizations of demand up and including that period.}

Then, our robust cost measure for selecting a reserve-crew schedule $\theta \in \Theta_n$ is given by 

\begin{align}
    c_\theta = (c^\textsc{s}\sum_{t \in \cal{T}}\xi_t(\theta) \ - \ c_0),
\end{align}
where $c^\textsc{s} $ is a fixed parameter denoting the costs of using reserve-schedule $\theta$, and $c_0$ is the robust costs corresponding to the current published schedule. For the first published schedule, we let $c_0 = 0$.

In Algorithm \ref{alg:evaluation}, we denote how to calculate the $\sum_{t\in\cal{T}} \xi_t$ for a given reserve schedule. Here, we calculate this fraction iteratively and denote its cumulative sum with $\Xi(\theta)$. We calculate $Q^\theta_t$ in a similar way, i.e., we calculate it iteratively and denote it with $Q(\theta)$.

\begin{algorithm}[t]
\footnotesize
  \caption{Evaluation of a reserve schedule}
  \label{alg:evaluation}
    \hspace*{\algorithmicindent} \textbf{Output:} Expected shortfall of reserve schedule $\Xi(\theta)$
  \begin{algorithmic}[1]
  \Procedure{Evaluate}{$\theta$, $\cal{S}$, $\cal{T}$, $f_t$}
    \State $\Xi(\theta) \gets 0$
    \State $Q(\theta) \gets (1, 0, \dots, 0)$
    \For{$t \in \cal{T}$}
        \State $\Xi(\theta) \gets \Xi(\theta) + \xi_{t}(\theta)$ 
        \State $Q(\theta) \gets Q(\theta)\cdot P^{C\theta}_{t}$
        \State $Q(\theta) \gets Q(\theta)\cdot P^{A\theta}_{t}$
        \State $Q(\theta) \gets Q(\theta)\cdot P^{B\theta}_{t}$
    \EndFor
    \State \textbf{return} $\Xi(\theta)$
    \EndProcedure
  \end{algorithmic}
\end{algorithm}

\section{The branch-and-price algorithm} \label{sec:pricingproblem}
In this section, we present the branch-and-price algorithm for solving the RCRP. First, we detail the branch-and-price algorithm by discussing branching rules.
Second, we introduce the pricing problem and its solution procedure.

\subsection{Branching rules}
We propose a sequential branching approach. First, we require all cancellations to be integral. Cancellations are typically highly costly and have a large effect on the resulting solution. By branching on these variables early, the number of flights which are to be covered becomes fixed and the resulting question is restricted to which crew member is responsible for taking on the flights to be covered. If multiple cancellation variables are fractional, branching is done on the most infeasible variable (i.e., the variable with value closest to $0.5$). Secondly, we force the deadheading variables to be integral. By doing so, we also restrict the pairing variables which include the flight leg to be deadheaded on to be integral. Again, branching is done on the most infeasible variable first. Finally, Ryan-Foster branching, as introduced by \citet{ryan1981}, is applied. 
In this branching rule, two successive flights in a fractional pairing (i.e., follow-ons) are forced to follow each other in a pairing on the one branch but not allowed to be in the same pairing on the other branch. As in the previous branching rule, follow-ons from the most infeasible pairing are selected to branch on first. Define the set of fractional pairing variables for crew member $k$ as $\bar P^{k}= \{p\in P^k \rvert x_p^k \notin \mathbb{Z}\}$. Then we branch on the follow-on connection $(i,j)$ that maximizes
\begin{equation}
    \min{\left(\sum_{k\in K}\sum_{p\in \bar P^k \rvert (i,j)\in p}x^k_p, 1 - \sum_{k\in K}\sum_{p\in \bar P^k \rvert (i,j)\in p}x^k_p\right)}.
\end{equation}
Besides branching, node selection strategies may play an important role in branch-and-price algorithms. Preliminary experiments have shown that  different branching rules do not show a big impact on the performance and the branch-and-bound node with the best-bound is not outperformed by other node selection strategies. We, therefore, choose this approach in our branch-and-price algorithm.  

\subsection{The pricing problem}
The pricing problem consists of finding crew pairings of negative reduced cost $\hat{c}^k_p$, for each crew member $k \in \cal{K}$. Note that the pricing problem and accompanying solution procedure are both valid for regular-crew members and reserve-crew members. In addition, note that every individual crew member has a specific set of duty-legality rules, as this depends on the overall crew schedule for multiple days. Therefore, it is crucial to solve the pricing problem for each crew-member individually. 

A crew pairing $ p \in \cal{P}^k$ is  described on a directed graph $\Graph^k = (\mathcal{V}^k, \mathcal{A}^k)$, where the vertex set $\mathcal{V}^k$ denotes airports at different moments in time, and the arc set $\mathcal{A}^k$ denotes the transitions between the airports (or the stay at the same airport). 

The vertex set $\mathcal{V}^k$ is further partitioned into a source node $v_s^k \in \cal{V}^k$ corresponding to the earliest possible starting time of crew member $k$, a sink node $v_e^k \in \cal{V}^k$ corresponding to the airports where crew member $k$  may end their day, 
and a set of assignment nodes $\mathcal{V}^\textsc{A} \subset \mathcal{V}^k$ denoting departure or arrival airports at a specific time so that flight legs can be represented by arcs between those.

\begin{figure}[t] 
\centering
\begin{tikzpicture}[->,>=stealth',shorten >=1pt,auto,
semithick, every text node part/.style={align=center}, scale = 0.9]
  \tikzstyle{every state}=[draw=black]
  \tiny
  \node[state]      (Source)                                       {AMS \\  $\ge 06:00$};
  
  \node[state]      (1_2)[above right = 0.5cm and 1cm of Source, text width=1cm] { RTM \\ 07:00};
   \node[state]      (1_1)[above = 0.5cm of 1_2, text width=1cm] { RTM \\ 06:55};
    \node[state]      (1_4)[below right = 0.5cm and 1cm of Source, text width=1cm] { AMS \\ 06:30};
    \node[state]      (1_5)[below = 0.5cm of 1_4, text width=1cm] { AMS \\ 06:35};
    \node[state]      (2_1)[right = 1cm of 1_1, text width=1cm] { BCN \\ 8:55};
    \node[state]      (2_2)[right = 1cm of 1_2, text width=1cm] { BCN \\ 9:00};
    \node[state]      (2_4)[right = 1cm of 1_4, text width=1cm] { ALC \\ 9:10};
    \node[state]      (2_5)[right = 1cm of 1_5, text width=1cm] { ALC \\ 9:15};
    \node[state]      (3_1)[right = 1cm of 2_1, text width=1cm] { BCN \\ 9:25};
    \node[state]      (3_2)[right = 1cm of 2_2, text width=1cm] { BCN \\ 9:30};
    \node[state]      (3_4)[right = 1cm of 2_4, text width=1cm] { ALC \\ 9:45};
    \node[state]      (3_5)[right = 1cm of 2_5, text width=1cm] { ALC \\ 9:50};
    \node[state]      (4_1)[right = 1cm of 3_1, text width=1cm] { AMS \\ 11:25};
    \node[state]      (4_2)[right = 1cm of 3_2, text width=1cm] { AMS \\ 11:30};
    \node[state]      (4_4)[right = 1cm of 3_4, text width=1cm] { AMS \\ 12:25};
    \node[state]      (4_5)[right = 1cm of 3_5, text width=1cm] { AMS \\ 12:30};
    \node[state]      (5_1)[right = 1cm of 4_1, text width=1cm] { AMS \\ 12:15};
    \node[state]      (5_2)[right = 1cm of 4_2, text width=1cm] { AMS \\ 12:20};
    \node[state]      (6_1)[right = 1cm of 5_1, text width=1cm] { GRQ \\ 13:00};
    \node[state]      (6_2)[right = 1cm of 5_2, text width=1cm] { GRQ \\ 13:05};
    \node[state]      (Sink)[right = 1.8cm of 4_4, text width=1.1cm] {AMS \\ $\le$ 15:30};

  \path(Source) edge       node[sloped,pos=0.8] {1:00}  (1_2)
    (Source) edge       node[sloped,pos=0.3] {0}  (1_4)
    (Source) edge       node[sloped,pos=0.2] {0}  (1_5)
    (1_1) edge       node[sloped,pos=0.5] {2:00}  (2_1)
    (1_2) edge       node[sloped,pos=0.5] {2:00}  (2_2)
    (1_4) edge       node[sloped,pos=0.5] {2:40}  (2_4)
    (1_5) edge       node[sloped,pos=0.5] {2:40}  (2_5)
    (2_1) edge       node[sloped,pos=0.5] {0:30}  (3_1)
    (2_2) edge       node[sloped,pos=0.5] {0:30}  (3_2)
    (2_4) edge       node[sloped,pos=0.5] {0:35}  (3_4)
    (2_5) edge       node[sloped,pos=0.5] {0:35}  (3_5)
    (2_5) edge       node[sloped,pos=0.5] {0:30}  (3_4)
     (2_4) edge       node[sloped,pos=0.05] {0:40}  (3_5)
     (2_1) edge       node[sloped,pos=0.2] {0:35}  (3_2)
     (3_1) edge       node[sloped,pos=0.5] {2:00}  (4_1)
    (3_2) edge       node[sloped,pos=0.5] {2:00}  (4_2)
    (3_4) edge       node[sloped,pos=0.5] {2:40}  (4_4)
    (3_5) edge       node[sloped,pos=0.5] {2:40}  (4_5)
    (4_4) edge       node[sloped,pos=0.4] {0}  (Sink)
    (4_5) edge       node[sloped,pos=0.5] {0}  (Sink)
    (4_1) edge[bend right]       node[sloped,pos=0.5] {0}  (Sink)
    (4_2) edge[bend right]       node[sloped,pos=0.5] {0}  (Sink)
    (4_1) edge       node[sloped,pos=0.5] {0:50}  (5_1)
    (4_2) edge       node[sloped,pos=0.5] {0:50}  (5_2)
    (4_1) edge       node[sloped,pos=0.05] {0:55}  (5_2)
     (4_2) edge       node[sloped,pos=0.5] {0:45}  (5_1)
     (5_1) edge       node[sloped,pos=0.5] {0:45}  (6_1)
     (5_2) edge       node[sloped,pos=0.5] {0:45}  (6_2)
     (6_1) edge [bend right]      node[sloped,pos=0.7] {2:30}  (Sink)
    ;

\end{tikzpicture}
\caption{Illustrative example of a pricing problem for crew member $k$ corresponding to the illustrative flight schedule in Table \ref{tab:exampleflights}}
\label{fig:graphstructure}

\end{figure}
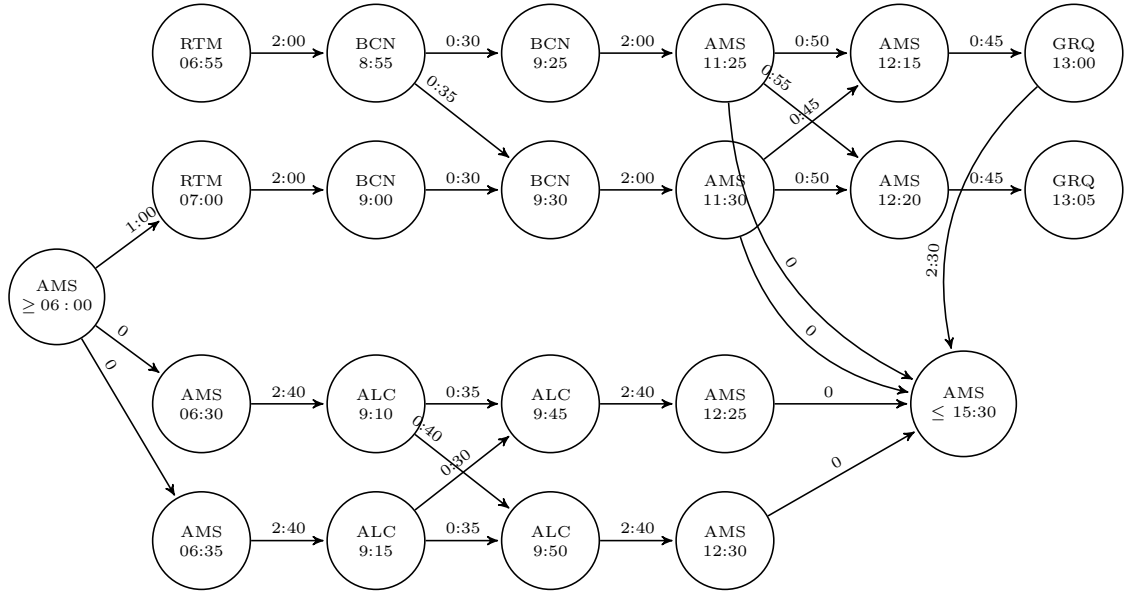

Finding a crew pairing of negative reduced cost equals solving a Resource Constrained Shortest Path Problem (RCSPP) in $\Graph^k$. Consider the following illustrative example, accompanying to Figure \ref{fig:graphstructure}, to clarify the structure of $\Graph^k$.

\begin{example}
In Figure \ref{fig:graphstructure}, we consider an example graph $\Graph^k$ for an arbitrary crew member $k \in \cal{K}$. Let the flight schedule be given as denoted in Table \ref{tab:exampleflights}. Let the shift of reserve-crew member $k$ start at 6:00 at AMS airport. Consider five flights as given in Table \ref{tab:exampleflights}. Transportation between the bases AMS, RTM and GRQ is possible by public transport and every flight is allowed to be re-timed such that they depart $5$ minutes earlier. In this example, crew member $k$ should arrive before 15:30 to not violate his minimum rest time. As a consequence, the arc between GRQ at 13:05 to AMS at 15:30 is infeasible (as the travel time is 2:30h) and therefore left out of the graph. The minimum sit time equals $30$ minutes, which renders some copy combinations infeasible. For instance, there is no arc between arriving at BCN at 9:00h and leaving BCN at 9:25, as this would only result in 25 minutes rest time between successive flight legs.  
\end{example}

\begin{table}[!b]
\centering
\caption{Example flight schedule}
\begin{tabular}{cccc}
    \FL DEP & ARR & DEP Time & ARR Time 
    \ML AMS & ALC & 6:35 & 9:15 
    \NN ALC & AMS & 9:50 & 12:30
    \NN RTM & BCN & 7:00 & 9:00
    \NN BCN & AMS & 9:30 & 11:30
    \NN AMS & GRQ & 12:20 & 13:05
    \LL 
\end{tabular}
\label{tab:exampleflights}
\end{table}

\subsection{Solving the pricing problem} \label{sec:pulse}
In order to solve RL-RCRP, we need to solve its corresponding pricing problems \eqref{eq:p1} and \eqref{eq:p2}, and add the generated crew pairings to RL-RCRP until no crew pairings of negative reduced cost can be found. Then, the corresponding dual solution is optimal, and we solved RL-RCRP to optimality in the current branch-and-bound node. Note that we enumerate the reserve-crew schedules up front, i.e., they do not play a part in the pricing problem. 

The pricing problem which we need to solve in the RCRP is a Resource Constrained Shortest Path Problem (RCSPP) that does not exhibit linear patterns in the use of resources. To not rely on dominance criteria within classical labeling algorithms, we propose a { pulse algorithm} \citep[see, e.g.,][]{lozano2015, schrotenboer2018} to solve the RCSPP. Besides that it does not rely on dominance criteria, it is simple but effective in its use. In summary, the pulse algorithm contains in a depth-first search through the graphs $\Graph^k$ while using pre-processed upper bounds in order to effectively prune the search. The procedures to solve the pricing problems are given in Algorithms \ref{alg:pricing} - \ref{alg:pulse}. {In Algorithm 2, the general procedure of solving our pricing problem is considered. It concerns a call to a bound procedure (i.e., Algorithm 3) and a call to the pulse algorithm (i.e., Algorithm 4). In the following, we discuss these algorithms one by one.}

The pulse algorithm {(Algorithm 4)} can be best characterized as a pulse propagating up and down in a depth-first search tree. At every node, the algorithm either continues the partial path $\Path$ it is following or is pruned. The resource consumption for a partial path $\Path$ is denoted using $\Phi_1(\Path)$ and $\Phi_2(\Path)$ for flying- and duty time, respectively. Reduced cost for a partial path $\Path$ is denoted using $r(\Path)$.

The pricing problem described in Algorithm \ref{alg:pricing} constitutes of two phases. In the first phase, a bounding scheme (refer to Algorithm \ref{alg:bound}) calculates lower bounds on the reduced costs.
A step size parameter $\Delta_p$ determines how many of these bounds are calculated. The bounding phase can be done independently of the subproblem for crew number $k$ (that is, the bounding phase holds for all crew members $k \in K$) and is described using a lower bound matrix $B = \left[b(v_i, \tau)\right]$, where $v_i$ corresponds to a given node and $\tau$ to the amount of resources spent. Secondly, the subproblem for crew member $k$ is solved using the pulse algorithm with a depth-first search of the solution space, where sub-optimal pairings are pruned using a combination of feasibility pruning (eliminating partial pairings which exceed the duty or flying time limits) and lower bound pruning using the lower bounds as calculated in the previous phase. Finally, if the pulse algorithm finds any pairing with negative reduced cost in subproblem $k\in K$, the pairing with the most negative reduced cost is added to the restricted master problem. That is, at each round of pricing, at most $k$ new pairings are added to the problem. 

\begin{algorithm}
\footnotesize
\caption{Pricing problem}
\label{alg:pricing}
\hspace*{\algorithmicindent} \textbf{Output:} New columns for the restricted master problem (RL-RCRP)
\begin{algorithmic}[1]
\Procedure{Pricing}{$\Delta_p$}
\State $\operatorname{bound}(\Delta_p, \phi_1^L)$
\For{$k \in K / K^R$}
    \State $\Path \gets \{\}$
    \State $r(\Path) \gets (-\delta^k + c^\textsc{a})$
    \State pulse$(v_s^k,0, r(\Path), 0, \Path)$
    \If{$r(\Path^*) < 0$}
        \State $\operatorname{addPairing}(\Path^*)$
    \EndIf 
\EndFor
\For{$k \in K^R$}
    \State $\Path \gets \{\}$
    \State $r(\Path) \gets (-\delta^k - \epsilon^k)$
    \State pulse$(v_s^k, 0, r(\Path), 0, \Path)$
    \If{$r(\Path^*) < 0$}
        \State $\operatorname{addPairing}(\Path^*)$
    \EndIf 
\EndFor
\EndProcedure
\end{algorithmic}
\end{algorithm}
 
We further elaborate on the pulse procedure {(i.e., Algorithm 4)}. If we intend to add node $v_i$ to the current partial path $\Path$, the function $\operatorname{isFeasible}$ checks whether adding this node does not exceed flight- or duty time limitations or violates feasibility in the monthly schedule. The function $\operatorname{checkBounds}$ checks whether going to this node within the current path is guaranteed to be sub-optimal given the current amount of flying spent $\Phi_1(\Path)$. That is, it prunes the search tree whenever $r(\Path) + b\left[v_i, \Phi_1(\Path)\right] \ge r(\Path^*),$ where $r(\Path^*)$ corresponds to the current best found solution. 

Bounds are only calculated for a finite number of possible amounts of flying time used, as determined by the parameter $\Delta_p$. In Algorithm \ref{alg:bound}, $f_i$ and $d_i$ denote the incurred flight and duty time when adding node $i$ to the current partial path. We ignore constraints on the duty time while constructing the lower bounds. In addition, note that the bound procedure is independent of the actual crew member $k \in \cal{K}$. Hence, when  $\tau = 0$, see Algorithm \ref{alg:bound}, the bound procedure does not provide us with optimal solutions to the RCSPP.

If no bound is known for $\Phi_1(\Path)$, one should round $\Phi_1(\Path)$ down to the highest amount of flying time which is in the bound matrix. When adding the current node $v_i$ is feasible and possibly optimal, we attempt to add a further node to the current partial path of the set $U^k_i = \{v_j \in \Edge^k | (v_i, v_j)\in \Arc^k\}$. Whenever a partial path $\Path$ reaches the sink node $v_e^k$ with minimal reduced cost so far, we set $\Path^* = \Path$ and $r(\Path^*) = r(\Path)$.

  \begin{algorithm}
  \footnotesize
  \caption{ Bound procedure}
  \label{alg:bound}
   \hspace*{\algorithmicindent} \textbf{Output:} Lower bound matrix $B = [b(v_i, \tau)]$
  \begin{algorithmic}[1]
  \Procedure{bound}{$\Delta_p$, $\phi_1^L$}
  \State $\tau \gets \phi_1^L$
    \While{$\tau \ge 0$}
        \State $\tau \gets \tau - \Delta_p$
        \State $\Path \gets \{\}$
        \For {$v_i \in \Edge^A$}
        \State $\operatorname{pulse}$($v_i$, $0$, $\tau$, $0, \Path$)
        \If {$\Path^* = \{\}$}
            \State $[b(v_i, \tau)] \gets \infty$
        \Else 
            \State $[b(v_i, \tau)] \gets r(\Path^*)$
        \EndIf
            
        \EndFor
     \EndWhile
    \EndProcedure
  \end{algorithmic}
  \end{algorithm}

\begin{algorithm}
\footnotesize
  \caption{Pulse procedure}
  \label{alg:pulse}
  \begin{algorithmic}[1]
  \Procedure{Pulse}{$v_i$, $r(\Path)$, $\Phi_1(\Path)$, $\Phi_2(\Path), \Path$}
    \If{$\operatorname{isFeasible}$($v_i$, $\Phi_1(\Path)$, $\Phi_2(\Path)$, $\Path$)}
        \If {not $\operatorname{checkBounds}$($v_i$, $\Phi_1(\Path)$, $r(\Path)$)}
        \State $\Path^{'} \gets \Path \cup \{v_i\}$
        \State $\Phi_1(\Path^{'}) \gets \Phi_1(\Path) + f_i$
        \State $\Phi_2(\Path^{'}) \gets \Phi_2(\Path) + d_i$
        \For{$v_j \in U^k_i$}
            \State $r(\Path^{'}) \gets r(\Path) + r_{ij}$
            \State $\operatorname{pulse}(v_j, r(\Path^{'}), \Phi_1(\Path^{'})$, $\Phi_2(\Path^{'}), \Path^{'})$
        \EndFor
       \EndIf
    \EndIf
    \EndProcedure
  \end{algorithmic}
  \end{algorithm}

\section{Case study} \label{sec:experimentaldesign}
The design of our model is based on a real-life case from a medium-sized Dutch international airline. We validate our model and solution technique based on the data of this airline. We follow the progress of their crew schedule over the tracking period and the disruptions which they encounter. We then let our model recover their schedules accordingly. We compare the results with a benchmark model, as we describe in Section \ref{sec:benchmark}.

\subsection{Benchmark model: Traditional Crew Recovery Problem} \label{sec:benchmark}
We compare our results with that of a Traditional Crew Recovery Problem (TCRP). The TCRP, unlike our model, uses no penalties for unfavorable characteristics to assess new pairings and penalizes reserve usage through a single cost parameter $c_s$ only. Hence, this model focuses on the strict minimization of planned costs only. For a fair comparison, the TCRP is restricted to using a maximum number of $n$ reserves when comparing its performance to that of the RCRP(n). 
The resulting MIP formulation can be found in Appendix \ref{sec:TCRPformulation}.
We show, through the upcoming experiments, that by applying a more reliable approach during the tracking period, results on the day of execution can be improved.  

\subsection{Experimental setup} \label{sec:data}
We consider two sets of each 10 instances. The first set of instances concerns 175 flights (39 aircraft) with 98 available crew-members, of which 19 are scheduled for a reserve-shift. The second set of instances concerns 309 flights (75 aircraft) with 169 crew-members, of which 29 are assigned to be reserve crew. The number of included bases (i.e., airports where the crew members start and end their day) differs between the experiments: it equals 4 for the first set and 1 for the second set of experiments. 

\begin{figure}[!b] 
\centering
\begin{tikzpicture}[->,>=stealth',shorten >=1pt, semithick, scale = 0.4]
\tiny
  \tikzstyle{every state}=[draw=black]

  \node[state,scale = 1, align=center]      (InputCrew)                                       {\textbf{Crew}\\ \textbf{Schedule}};
  \node[state,scale = 1, align=center]      (InputFlight)[below = 1cm of InputCrew] {\textbf{Flight} \\ \textbf{Schedule}};
  \node[shape = rectangle, draw, align=center, anchor = east, scale = 1, minimum width=1cm,minimum height=1cm, text width = 3cm]      (Disruptions)[below right = 1.38cm and 1cm of InputCrew, text width=2cm] {\textbf{Disruptions}};
   \node[draw, scale = 1, shape = diamond, align = center, minimum width=2cm,minimum height=2cm] (Tracking)[below right = 4cm and 1.5cm of InputCrew]{\textbf{Tracking} \\
   \textbf{period}?};
   \node[shape = rectangle, draw, align=center, anchor = east, scale = 1, minimum width=1cm,minimum height=1cm, text width = 3cm]      (SolveTCRP)[right = 1cm of Tracking, text width=2cm] {\textbf{Solve} \\ \textbf{TCRP}};
   \node[shape = rectangle, draw, align=center, anchor = east, scale = 1, minimum width=1cm,minimum height=1cm, text width = 3cm]      (SolveRCRP)[left = 3.5cm of Tracking, text width=2cm] {\textbf{Solve} \\ \textbf{RCRP}};
   \node[state,scale = 1, align=center]      (OutputCrew)[left = 1cm of InputCrew]                                       {\textbf{Updated} \\ \textbf{Crew}\\ \textbf{Schedule}};
  \node[state,scale = 1, align=center]      (OutputFlight)[left = 1cm of InputFlight] {\textbf{Updated} \\ \textbf{Flight} \\ \textbf{Schedule}};
  \node[fill = gray, shape = rectangle, draw, align=center, anchor = east, scale = 1, minimum width=1cm,minimum height=1cm, text width = 3cm]      (EndExperiment)[above = 1.6cm  of SolveTCRP, text width=2cm] {\textbf{End}};

  \path(InputCrew) edge       node[sloped,pos=1] {}  (Disruptions)
  (InputFlight) edge       node[sloped,pos=1] {}  (Disruptions)
  (Disruptions) edge       node[sloped,pos=1] {}  (Tracking)
  (Tracking) edge       node[sloped,pos=0.5, above] {\textbf{Yes}}  (SolveRCRP)
  (Tracking) edge       node[sloped,pos=0.5, above] {\textbf{No}}  (SolveTCRP)
  (SolveRCRP) edge  [bend left=20] node[sloped,pos=0.6, above] {}  (OutputCrew)
  (SolveRCRP) edge       node[sloped,pos=0.3, above] {}  (OutputFlight)
 (OutputCrew) edge      node[->,sloped,pos=0.5, above] {}  (InputCrew)
 (OutputFlight) edge   node[->,sloped,pos=0.5, above] {}  (InputFlight)
  (SolveTCRP) edge       node[->,sloped,pos=0.5, above] {}  (EndExperiment)

    ;

\end{tikzpicture}
\caption{Illustration of the experimental set-up}
\label{fig:experimentaldesign} 
\end{figure}
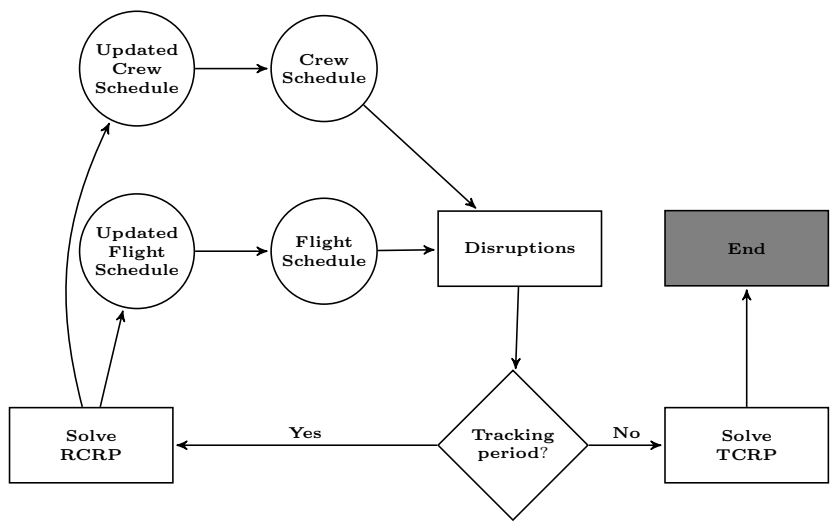

In Figure \ref{fig:experimentaldesign}, we illustrate our experimental design. We consider a tracking period of 30 days throughout which disruptions might occur. In the first experimental setup, as a test case, we consider all the disruptions throughout the 30 days and only update the schedule once. In the second experimental setup, we partition the tracking period in four periods of equal length and update the schedule four times (corresponding to a weekly update).
Hereby, the first experiment serves as a validation of the RCRP in relation to the TCRP, i.e., the differences in actual results should be comparable though the RCRP should perform more reliable. The second set of experiments is in line with practice, where crew schedules are updated in a weekly basis in order to inform crew members timely about the disruptions in their schedule.

The probability of a single disruption (for a flight leg not being covered) is denoted by {}{$\psi_i^j$} for each instance $i$, where $j$ denotes the period in which we recover from the accumulated disruptions. For example, $\psi^{0.25T}_i$ means the probability of a single disruption occurring in the first quarter of the tracking period. For the first experimental setup, we only consider $\psi_i^T$, as we accumulate all disruptions throughout the tracking period. We denote the probability of a disrupted flight leg on the day of execution with $\psi_i^O$. This implies that, for 
a period $j$ with $b_j$ flights departing from the home base(s), the demand for reserve personnel
can be described using a binomial distribution $B(b_j;\psi_i^j)$. In Table \ref{tab:scenariodisruption}, we provide the number of resulting disruptions, denoted by $\kappa_i^T$ and $\kappa_i^O$, for each recovering moment during the tracking phase. We named the instances according to the experimental setup they belong to ($E1_i$ or $E2_i$). In summary, instance set E1(1-10) is primarily used to demonstrate the effect of using the RCRP for different set of scenarios based on probabilities $\psi_i^T$ and $\psi_i^O$. This provides insights in how sensitive the model is to disruptions happening at different configurations. In instance set E2(1-10), the probability of disruptions over time is based on historical data. Furthermore, it more closely resembles a real life situation in which disruptions are periodically monitored. 

\begin{table}[htbp]
\footnotesize
\caption{Scenario descriptions and amount of disrupted (reserve) personnel}
\centering
\begin{tabular}{l|rrrrrrrrrr}
    \FL Instance & $\psi^{0.25T}_i$ & $\psi^{0.5T}_i$ & $\psi^{0.75T}_i$ & $\psi^T_i$ & $\psi^O_i$ &  $\kappa^{0.25T}_i$ & $\kappa^{0.5T}_i$ & $\kappa^{0.75T}_i$ & $\kappa^T_i$ & $\kappa^O_i$ 
    \ML $E1_1$  & - & - & - & 0.10 & 0.10  & - & - & -& 14 & 6 
    \NN $E1_2$  & - & - & - & 0.10 & 0.10  & - & - & -& 8 & 11
    \NN $E1_3$  & - & - & - & 0.10 & 0.08  & - & - & -& 8 & 9
    \NN $E1_4$ & - & - & - & 0.10 & 0.08  & - & - & -& 12 & 5
    \NN $E1_5$  & - & - & - & 0.08 & 0.10  & - & - & -& 2 & 12
    \NN $E1_6$  & - & - & - & 0.08 & 0.10  & - & - & -& 7 & 7
    \NN $E1_7$ & - & - & - & 0.05 & 0.15 & - & - & - & 5 & 11
    \NN $E1_8$ & - & - & -  & 0.05 & 0.15 & - & - & - & 4 & 15
    \NN $E1_9$ & - & - & -  & 0.15 & 0.05  & - & - & -& 18 & 4
    \NN $E1_{10}$  & - & - & -  & 0.15 & 0.05  & - & - & -& 20 & 7 \\[0.1cm]
    $E2_1$ & 0.02 & 0.01 & 0.01 & 0.05 & 0.06 & 3 & 1  & 2 & 5 & 14
	\NN $E2_2$ & 0.02 & 0.01 & 0.01 & 0.05 & 0.06& 5 & 4 & 1 & 6 & 12
	\NN $E2_3$ & 0.02 & 0.01 & 0.01 & 0.05 & 0.06& 5 & 0 & 0 & 5 & 12
	\NN $E2_4$ & 0.02 & 0.01 & 0.01 & 0.05 & 0.06& 4 & 1 & 2 & 7 & 9
	\NN $E2_5$ & 0.02 & 0.01 & 0.01 & 0.05 & 0.06& 2 & 3 & 0 & 8 & 12
	\NN $E2_6$ & 0.02 & 0.01 & 0.01 & 0.05 & 0.06& 3 & 3 & 2 & 10 & 12
	\NN $E2_7$ & 0.02 & 0.01 & 0.01 & 0.05 & 0.06& 3 & 2 & 2 & 10 & 12
	\NN $E2_8$ & 0.02 & 0.01 & 0.01 & 0.05 & 0.06& 4 & 0 & 0 & 6 & 12
	\NN $E2_9$ & 0.02 & 0.01 & 0.01 & 0.05 & 0.06& 4 & 1 & 1 & 9 & 12
	\NN $E2_{10}$ & 0.02 & 0.01 & 0.01 & 0.05 & 0.06 & 6 & 2 & 4 & 4 & 8
	\LL
\end{tabular}
\label{tab:scenariodisruption}
\end{table}

All published duties are subject to change, but corresponding to labor regulations, the starting time of a regular crew member's duty may be at most two hours earlier than its originally planned duty. Furthermore, we set $\phi_1^L = 9:00h$, $\phi_1^S = 12:00h$, $\phi_2^L = 11:00h$, $\phi_2^S = 13:00h$, $\phi_3 = 0:30h$, and $\phi_4 = 10:00h$. Recall that these parameters model the duty legality rules, as explained in Section \ref{sec:propmodel}. All the values of the cost parameters can be found in Table \ref{tab:costsparameters}. 

Finally, we consider 5 flight copies {}{(with intervals of 5 minutes)} for each flight leg to model the re-timings. In addition, we limit the number of reserve schedules to 
 $$ \max_{m} \{m | \sum_{i = 0}^m \frac{M!}{i!(M-i)!}\le \Pi\},$$ where {}{$M$ equals the total number of reserves} and $\Pi$ is a parameter limiting the amount of reserve schedules to be evaluated. We set $\Pi = 100,000$ for Experiment 1 and to $500,000$ for Experiment 2. 
The pricing problem is solved with a pulse parameter of $\Delta_p = 60$.

\begin{table}[htbp]
	\centering
	\footnotesize

	\caption{Costs parameters relative to a minute of flying time}
	\begin{tabular}{lll}
		\FL Parameter & Value & Explanation 
		\ML $c^{d}_{\omega}$ & \{250, -1750\} & Deadheading costs for original (250) and re-timed flight copies (-1750)
		\NN $c_f^\textsc{c}$ & 1000000 & Flight cancellation costs
		\NN $c^\textsc{s}$ & \{2500, 360\} & Per unit cost of expected shortfall (RCRP - 2500) and reserve-crew cost (TCRP - 360)
		\NN $c^\textsc{r}_f$ & 2000 & Re-timing costs
		\NN $\xi$ & 5/8 & Parameter to compare duty and flight time
		\NN $PC^{\min}$ & 360 & Minimum pay-and-credit for a crew member
		\NN $c_\textsc{t}$ & 60 & Fixed transportation costs between bases
		\NN $c_\textsc{a}$ & 200 & Schedule altering costs
		\NN $g_1$ & 360 & Maximum penalty for flying time regulations
		\NN $\Delta_1$ & 4 & Per unit cost for exceeding maximum flying time
		\NN $g_2$ & 360 & Maximum penalty for duty time regulations 
		\NN $\Delta_2$ & 4 & Per unit cost for exceeding maximum duty time
		\NN $g_3$ & 180 & Maximum penalty for rest time regulations
		\NN $\Delta_3$ & 2 & Per unit cost for exceeding minimum rest time
		\NN $g_4$ & 120 & Maximum penalty for sit time regulations
		\NN $\Delta_4$ & 12 & Per unit cost for exceeding minimum sit time
		\LL 
	\end{tabular}
	\label{tab:costsparameters}
\end{table}

\subsection{Results} \label{sec:results}
In this section, we present the results of solving the instances of both Experiment 1 and Experiment 2. Each instance is solved as indicated in Figure \ref{fig:experimentaldesign}, and this solution is referred to as the RCRP solution. For benchmark purposes, we calculated the solutions without explicitly considering reliable reserve-crew schedules, and we call these TCRP solutions.  For each instance, we provide the statistics of two solutions. First we focus on the proposed solutions at the end of the tracking period, and second, we present the solutions after recovering from the disruptions on the day of execution. 

The branch-and-price algorithm is programmed in C++ by using SCIP 4.0 \citep{scip} in combination with SoPlex 3.0.1 to calculate the LP relaxations. All computational tests are performed on a Linux virtual machine running on a laptop with an Intel Core i5-7200 CPU with 3GB of random access memory.


At our case airline, the schedule is continuously maintained over the “Tracking” period (starting from the last 30 days before the day of execution). In the last 24 hours before execution, the maintenance of the schedule moves towards the `Control' phase, in which a different set of planners make sure the flights all go on as planned. In short, the focus in the tracking period is on “Maintaining”, while the focus on the day of operation generally lies on `Making sure we fly as much as possible'.  
The resulting crew schedule at the end of a tracking period (with characteristics defined in Table 4) thus form the input schedule for the day of operation (and as such impact the robustness of the schedule). During the day of operation, planners deal with additional disruptions occurring on the crew schedule. The resulting characteristics of dealing with these additional disruptions on the operation day are presented in Table 5. 
In Tables 4 and 5, the column headings indicate the number of alterations (Alt),  re-timings (Ret),  reserves used (Res), cancellations (Canc), aggregated penalties (Pen), expected shortfall (Shortfall), aggregated costs (costs) and solution times in seconds (Time).

\begin{table}[htbp]
\centering
\footnotesize
\caption{Solution characteristics at the end of the tracking period}
\begin{tabular}{l|r|lr|lr|lr|lr|lr|lr|lr|l}
    \FL  \multicolumn{1}{c}{} & \multicolumn{8}{r|}{RCRP } &  \multicolumn{8}{l}{TCRP}
    \ML Instance  & \multicolumn{2}{c}{Alt} & \multicolumn{2}{c}{Ret} & \multicolumn{2}{c}{Res} & \multicolumn{2}{c}{Canc} & \multicolumn{2}{c}{Pen} & \multicolumn{2}{c}{Shortfall} & \multicolumn{2}{c}{Costs}  & \multicolumn{2}{c}{Time} \\ \midrule
$E1_1$    & 9            & 7            & 0            & 0            & 9            & 9            & 0             & 0             & 3740          & 4620          & 2.01          & 2.81          & 31913          & 31953          & 11          & 2          \\
$E1_2$    & 0            & 0            & 0            & 0            & 6            & 6            & 0             & 0             & 1420          & 1380          & 0.79          & 1.01          & 32698          & 32578          & 8           & 0          \\
$E1_3$    & 0            & 0            & 0            & 0            & 6            & 6            & 0             & 0             & 1450          & 1670          & 0.91          & 1.31          & 32578          & 32578          & 8           & 1          \\
$E1_4$    & 3            & 1            & 0            & 0            & 9            & 9            & 0             & 0             & 2220          & 2450          & 2.01          & 3.19          & 32518          & 32383          & 10          & 2          \\
$E1_5$    & 0            & 0            & 0            & 0            & 2            & 2            & 0             & 0             & 1450          & 1610          & 0.29          & 0.29          & 32578          & 32578          & 5           & 0          \\
$E1_6$    & 2            & 0            & 0            & 0            & 7            & 7            & 0             & 0             & 1550          & 1980          & 2.34          & 3.29          & 32698          & 32818          & 8           & 2          \\
$E1_7$    & 1            & 0            & 0            & 0            & 4            & 4            & 0             & 0             & 1450          & 1450          & 0.40          & 0.87          & 32698          & 32578          & 4           & 0          \\
$E1_8$    & 1            & 0            & 0            & 0            & 4            & 4            & 0             & 0             & 1450          & 1450          & 0.39          & 0.68          & 32698          & 32578          & 6           & 0          \\
$E1_9$    & 11           & 11           & 1            & 1            & 9            & 9            & 4             & 4             & 5100          & 5120          & 2.06          & 4.11          & 30833          & 30543          & 7           & 3          \\
$E1_{10}$ & 11           & 11           & 3            & 3            & 9            & 9            & 9             & 9             & 5310          & 5430          & 4.10          & 5.01          & 29650          & 29854          & 9           & 2          \\[0.1cm]
$E2_1$    & 3            & 1            & 0            & 0            & 8            & 8            & 0             & 0             & 2270          & 2270          & 0.42          & 1.87          & 57775          & 57775          & 45          & 2          \\
$E2_2$    & 8            & 0            & 0            & 0            & 12           & 12           & 0             & 0             & 2330          & 2320          & 0.72          & 4.99          & 57775          & 57775          & 108         & 18         \\
$E2_3$    & 5            & 0            & 0            & 0            & 8            & 8            & 0             & 0             & 2240          & 2300          & 0.32          & 1.66          & 57835          & 57835          & 83          & 3          \\
$E2_4$    & 5            & 0            & 0            & 0            & 11           & 11           & 0             & 0             & 2300          & 2480          & 0.38          & 1.91          & 57895          & 57835          & 60          & 11         \\
$E2_5$    & 4            & 0            & 0            & 0            & 11           & 10           & 0             & 0             & 2450          & 2280          & 0.36          & 2.13          & 57655          & 57655          & 70          & 9          \\
$E2_6$    & 4            & 0            & 0            & 0            & 11           & 11           & 0             & 0             & 2450          & 2470          & 0.36          & 3.20          & 57655          & 57775          & 52          & 6          \\
$E2_7$    & 11           & 1            & 0            & 0            & 14           & 15           & 0             & 0             & 4000          & 3840          & 0.38          & 1.60          & 57149          & 57610          & 59          & 7          \\
$E2_8$    & 1            & 0            & 0            & 0            & 8            & 8            & 0             & 0             & 2140          & 2450          & 0.29          & 1.22          & 57835          & 57835          & 26          & 10         \\
$E2_9$    & 5            & 0            & 0            & 0            & 13           & 12           & 0             & 0             & 2520          & 2870          & 0.38          & 0.87          & 57415          & 57480          & 87          & 11         \\
$E2_{10}$ & 14           & 2            & 0            & 0            & 15           & 14           & 0             & 0             & 15188         & 4730          & 0.46          & 2.67          & 59594          & 57744          & 75          & 22         \\ \midrule
    \textbf{Avg. }     & \textbf{4.90} & \textbf{1.70} & \textbf{0.20} & \textbf{0.20} & \textbf{8.80} & \textbf{8.70} & \textbf{0.65} & \textbf{0.65} & \textbf{3151} & \textbf{2759} & \textbf{0.97} & \textbf{2.23} & \textbf{44972} & \textbf{44888} & \textbf{37} & \textbf{6} \\ \bottomrule
\end{tabular}
\label{tab:scenarioresultstracking}
\end{table}

In Table \ref{tab:scenarioresultstracking}, we observe that using the RCRP instead of the TCRP reduces the expected shortfall with 56.5\%, while the cost-increase is only 0.18\%. The cause can be found in the number of alterations, which are on average 4.9 and 1.7 for the RCRP and TCRP, respectively. Hence, additional alterations are taken in the RCRP in order to significantly lower the expected shortfall. 

In addition, we observe that when possible, both recovery models attempt to primarily recover the disrupted schedule through the use of reserve crew. Re-timing of flights is only considered as a last minute resort if cancellations cannot be prevented otherwise. For the TCRP, this also holds for altering existing schedules, if reserves cannot be plugged in directly to replace the disrupted crew member. For the RCRP, altering some (additional) schedules is a viable option for reducing future risk. The main reason for doing so arises through the objective of minimizing reserve schedule shortfall, which is significantly lower when using the RCRP than using the TCRP. Through altering a few extra schedules, the RCRP can use the reserves which contribute the least to schedule shortfall. Overall, the RCRP provides solutions with less aggregated penalties but it may actually allow additional penalties, if this would benefit the accompanying reserve schedule expected shortfall. 

In Table \ref{tab:scenarioresultsoperation}, we provide the solution characteristics at the end of the day of execution. Note that these are the results of a final round of disruptions on the solutions presented in Table \ref{tab:scenarioresultstracking} and afterward solving it with the TCRP.

From Table \ref{tab:scenarioresultsoperation}, it is observed that the average number of (last-minute) alterations is reduced from on average 4.35 to 2.75 if one uses the RCRP during the tracking period instead of the TCRP. It is also noted that the total penalty costs for the RCRP are lower than that of the TCRP, indicating that incorporating the robust measures during the tracking period (by solving the RCRP) results in more flexibility to avoid penalty costs on the day of execution.

\begin{table}[htbp]
\footnotesize

\centering
\caption{Solution characteristics at the end of day of execution}
\begin{tabular}{l|r|lr|lr|lr|lr|lr|lr|l}
    \FL  \multicolumn{1}{c}{} & \multicolumn{7}{r|}{RCRP } &  \multicolumn{7}{l}{TCRP}
    \ML Instance & \multicolumn{2}{c}{Alt} & \multicolumn{2}{c}{Ret} & \multicolumn{2}{c}{Res} & \multicolumn{2}{c}{Canc} & \multicolumn{2}{c}{Pen} & \multicolumn{2}{c}{Costs} & \multicolumn{2}{c}{Time} \\ \midrule
  $E1_1$    & 1             & 2             & 0            & 0             & 6             & 6             & 0            & 0             & 3980          & 3960          & 32643          & 32438          & 1           & 1           \\
$E1_2$    & 7             & 6             & 1            & 1             & 10            & 11            & 0            & 0             & 2490          & 2660          & 32927          & 32945          & 2           & 1           \\
$E1_3$    & 2             & 3             & 0            & 0             & 8             & 9             & 1            & 1             & 1760          & 1940          & 33000          & 33000          & 0           & 0           \\
$E1_4$    & 1             & 2             & 0            & 0             & 4             & 4             & 0            & 0             & 2250          & 2630          & 32758          & 32623          & 1           & 2           \\
$E1_5$    & 0             & 0             & 0            & 0             & 9             & 9             & 0            & 0             & 1690          & 1670          & 32938          & 32938          & 1           & 1           \\
$E1_6$    & 5             & 6             & 1            & 1             & 6             & 6             & 2            & 2             & 2770          & 3190          & 32510          & 32510          & 1           & 1           \\
$E1_7$    & 1             & 1             & 0            & 0             & 9             & 9             & 0            & 0             & 2340          & 2520          & 32638          & 32638          & 1           & 1           \\
$E1_8$    & 5             & 6             & 2            & 2             & 13            & 13            & 1            & 1             & 2260          & 2230          & 33321          & 33201          & 1           & 1           \\
$E1_9$    & 3             & 5             & 0            & 1             & 6             & 6             & 0            & 1             & 5390          & 5660          & 32045          & 32068          & 2           & 2           \\
$E1_{10}$ & 4             & 2             & 0            & 1             & 7             & 5             & 8            & 10            & 5830          & 5990              & 30190          & 29854               & 1           & 1           \\[0.1cm]
$E2_1$    & 9             & 12            & 4            & 3             & 11            & 9             & 2            & 4             & 5050          & 6650          & 57788          & 57049          & 131         & 161         \\
$E2_2$    & 0             & 6             & 0            & 0             & 12            & 9             & 0            & 0             & 2910          & 4760          & 57895          & 57629          & 21          & 37          \\
$E2_3$    & 2             & 5             & 0            & 0             & 10            & 9             & 0            & 0             & 2380          & 4760          & 57775          & 57475          & 37          & 48          \\
$E2_4$    & 3             & 5             & 0            & 0             & 8             & 7             & 0            & 0             & 3340          & 4710          & 57605          & 57470          & 28          & 60          \\
$E2_5$    & 3             & 7             & 0            & 0             & 11            & 10            & 0            & 0             & 4440          & 5140          & 57639          & 57564          & 54          & 145         \\
$E2_6$    & 3             & 5             & 0            & 0             & 11            & 8             & 0            & 0             & 4440          & 4680          & 57639          & 57389          & 59          & 63          \\
$E2_7$    & 0             & 1             & 0            & 0             & 10            & 8             & 0            & 0             & 4250          & 5090          & 57119          & 57480          & 24          & 46          \\
$E2_8$    & 2             & 3             & 0            & 0             & 10            & 10            & 0            & 0             & 3700          & 4210          & 57765          & 57765          & 19          & 25          \\
$E2_9$    & 0             & 1             & 0            & 0             & 10            & 10            & 0            & 0             & 3520          & 4050          & 57319          & 57355          & 19          & 53          \\
$E2_{10}$ & 4             & 9             & 0            & 0             & 6             & 7             & 0            & 0             & 13412         & 7540          & 58814          & 57440          & 16          & 81          \\ \midrule
   Avg.       & \textbf{2.75} & \textbf{4.35} & \textbf{0.40} & \textbf{0.45} & \textbf{8.85} & \textbf{8.25} & \textbf{0.70} & \textbf{0.95} & \textbf{3910} & \textbf{4108} & \textbf{45116} & \textbf{45736} & \textbf{21} & \textbf{37} \\ \bottomrule
   \end{tabular}
\label{tab:scenarioresultsoperation}
\end{table}

While one may argue that the reduction in alterations is due to more alterations during the tracking phase (see Table \ref{tab:scenarioresultstracking}), we argue that alterations on the final day are much more disruptive to the airline operations, creating delays and crew confusion. In addition, alterations close to execution have a strong negative impact on the crew-members' personal lives and their satisfaction.

Finally, reliable recovery during the tracking phase will even prevent cancellations during the final day of execution. In Instance $E2_1$, a grand total of $9$ reserves are used in the final solution under TCRP, but $4$ flights had to be canceled regardless. While all reserves are scheduled for a shift early in the morning, the disrupted flights are scheduled for departure late in the evening. In the schedule as produced by the RCRP, reserves are more spread out and as a result the cancellations are reduced. 
In summary, the advantage of the RCRP is the use of completely enumerated reserve-crew schedule, instead of just using reserve crew members without considering the impact future recovering on the flight execution day. Thus recovering to feasible flight schedules using the RCRP leads to less extreme events on flight execution day. 

As the instances of Experiment 2 are centered around a single basis, where crew starts and end their days, we conclude that reliable reserve-crew management is especially important for a carrier operating in a hub-and-spoke network, in which reserve crew members are located at the main hub from which the aircraft rotations are departing. Here, a significant number of last-minute alterations and even flight cancellations are prevented when one adopts the reliable reserve-crew scheduling approach. 

{}{We also test the computational performance of the our algorithm on Experiment 2 by increasing the number of flight copies. We report on the cost, total computation time, and the solution characteristics during the tracking phase and on the execution day in Tables \ref{tab:comp1} - \ref{tab:comp6}. As expected, with an increasing number of flight copies, the computation time of our algorithm increases. However, even for larger sizes of instances we observe that the RCRP provides better performance in terms of the cost efficiency than the TCRP. The computation time 11 flight copies stays reasonable, within 2h during the tracking phase and at most 1h during the execution day approximately.}

\begin{table}[htbp]
\centering
\footnotesize
\caption{Solution characteristics at the end of the tracking period for 7 flight copies}
\begin{tabular}{l|r|lr|lr|lr|lr|lr|lr|lr|l}
    \FL  \multicolumn{1}{c}{} & \multicolumn{8}{r|}{RCRP } &  \multicolumn{8}{l}{TCRP}
    \ML Instance  & \multicolumn{2}{c}{Alt} & \multicolumn{2}{c}{Ret} & \multicolumn{2}{c}{Res} & \multicolumn{2}{c}{Canc} & \multicolumn{2}{c}{Pen} & \multicolumn{2}{c}{Shortfall} & \multicolumn{2}{c}{Costs}  & \multicolumn{2}{c}{Time} \\ \midrule
$E2_1$    & 1           & 1   & 0       & 0   & 4             & 4   & 0         & 0 & 2370    & 2270 & 0.42      & 1.17  & 57775 & 57775 & 565  & 31   \\
$E2_2$    & 3           & 1   & 0       & 0   & 6             & 5   & 0         & 0 & 2310    & 3200 & 0.72      & 4.05  & 57775 & 57624 & 1014 & 524  \\
$E2_3$    & 1           & 0   & 0       & 0   & 3             & 4   & 0         & 0 & 2240    & 2300 & 0.34      & 1.66  & 57835 & 57835 & 1007 & 66   \\
$E2_4$    & 1           & 0   & 0       & 0   & 5             & 5   & 0         & 0 & 2420    & 2480 & 0.39      & 1.92  & 57895 & 57835 & 762  & 474  \\
$E2_5$    & 8           & 1   & 0       & 0   & 6             & 7   & 0         & 0 & 3760    & 2290 & 0.26      & 2.14  & 57364 & 57655 & 698  & 125  \\
$E2_6$    & 1           & 0   & 0       & 0   & 5             & 6   & 0         & 0 & 2670    & 2610 & 1.11      & 3.34  & 57775 & 57775 & 705  & 113  \\
$E2_7$    & 4           & 1   & 0       & 0   & 8             & 8   & 0         & 0 & 3100    & 3590 & 0.56      & 1.60  & 57535 & 57610 & 507  & 144  \\
$E2_8$    & 1           & 0   & 0       & 0   & 5             & 5   & 0         & 0 & 2140    & 2320 & 0.29      & 1.22  & 57835 & 57835 & 361  & 196  \\
$E2_9$    & 4           & 2   & 0       & 0   & 7             & 7   & 0         & 0 & 2520    & 2950 & 0.38      & 0.87  & 57415 & 57480 & 1020 & 175  \\
$E2_{10}$ & 2           & 0   & 0       & 0   & 4             & 4   & 0         & 0 & 4190    & 5930 & 0.39      & 1.35  & 57455 & 57354 & 1203 & 746  \\ \midrule
Avg.          & 2.6         & 0.6 & 0       & 0   & 5.3           & 5.5 & 0         & 0 & 2772    & 2994 & 0.49      & 1.93  & 57666 & 57678 & 784  & 259  \\ \bottomrule
\end{tabular}
\label{tab:comp1}
\end{table}

\begin{table}[htbp]
\centering
\footnotesize
\caption{Solution characteristics at the end of the tracking period for 9 flight copies}
\begin{tabular}{l|r|lr|lr|lr|lr|lr|lr|lr|l}
    \FL  \multicolumn{1}{c}{} & \multicolumn{8}{r|}{RCRP } &  \multicolumn{8}{l}{TCRP}
    \ML Instance  & \multicolumn{2}{c}{Alt} & \multicolumn{2}{c}{Ret} & \multicolumn{2}{c}{Res} & \multicolumn{2}{c}{Canc} & \multicolumn{2}{c}{Pen} & \multicolumn{2}{c}{Shortfall} & \multicolumn{2}{c}{Costs}  & \multicolumn{2}{c}{Time} \\ \midrule
$E2_1$    & 2           & 1   & 0       & 0   & 4             & 4   & 0         & 0 & 2270    & 2440 & 0.31      & 1.17  & 57775 & 57775 & 923  & 77   \\
$E2_2$    & 3           & 1   & 0       & 0   & 6             & 5   & 0         & 0 & 2310    & 3330 & 0.72      & 4.05  & 57775 & 57624 & 2569 & 971  \\
$E2_3$    & 1           & 0   & 0       & 0   & 3             & 3   & 0         & 0 & 2240    & 2140 & 0.30      & 1.86  & 57835 & 57835 & 2627 & 149  \\
$E2_4$    & 1           & 0   & 0       & 0   & 5             & 5   & 0         & 0 & 2420    & 2440 & 0.38      & 1.91  & 57895 & 57835 & 2147 & 1515 \\
$E2_5$    & 1           & 1   & 0       & 0   & 7             & 7   & 0         & 0 & 2340    & 2310 & 0.36      & 2.21  & 57655 & 57655 & 1552 & 383  \\
$E2_6$    & 1           & 0   & 0       & 0   & 5             & 6   & 0         & 0 & 2630    & 2470 & 1.11      & 3.34  & 57775 & 57775 & 1654 & 173  \\
$E2_7$    & 4           & 1   & 0       & 0   & 8             & 8   & 0         & 0 & 3100    & 3840 & 0.56      & 1.51  & 57535 & 57610 & 1216 & 514  \\
$E2_8$    & 1           & 0   & 0       & 0   & 5             & 5   & 0         & 0 & 2140    & 2480 & 0.29      & 0.68  & 57835 & 57835 & 888  & 531  \\
$E2_9$    & 3           & 2   & 0       & 0   & 7             & 7   & 0         & 0 & 2630    & 2850 & 0.38      & 0.82  & 57470 & 57480 & 2210 & 542  \\
$E2_{10}$ & 2           & 0   & 0       & 0   & 6             & 3   & 0         & 0 & 3800    & 6380 & 0.47      & 1.61  & 57705 & 57404 & 2149 & 2175 \\ \midrule
Avg.          & 1.9         & 0.6 & 0       & 0   & 5.6           & 5.3 & 0         & 0 & 2588    & 3068 & 0.49      & 1.92  & 57726 & 57683 & 1794 & 703  \\ \bottomrule 
          \end{tabular} 
\end{table}

\begin{table}[htbp]
\centering
\footnotesize
\caption{Solution characteristics at the end of the tracking period for 11 flight copies}
\begin{tabular}{l|r|lr|lr|lr|lr|lr|lr|lr|l}
    \FL  \multicolumn{1}{c}{} & \multicolumn{8}{r|}{RCRP } &  \multicolumn{8}{l}{TCRP}
    \ML Instance  & \multicolumn{2}{c}{Alt} & \multicolumn{2}{c}{Ret} & \multicolumn{2}{c}{Res} & \multicolumn{2}{c}{Canc} & \multicolumn{2}{c}{Pen} & \multicolumn{2}{c}{Shortfall} & \multicolumn{2}{c}{Costs}  & \multicolumn{2}{c}{Time} \\ \midrule
$E2_1$    & 2           & 1   & 0       & 0   & 4             & 4   & 0         & 0 & 2270    & 2270 & 0.31      & 1.96  & 57775 & 57775 & 2473 & 184  \\
$E2_2$    & 3           & 1   & 0       & 0   & 6             & 5   & 0         & 0 & 2310    & 3380 & 0.72      & 4.05  & 57775 & 57624 & 4396 & 4882 \\
$E2_3$    & 1           & 0   & 0       & 0   & 3             & 4   & 0         & 0 & 2240    & 2300 & 0.34      & 1.66  & 57835 & 57835 & 5850 & 1171 \\
$E2_4$    & 1           & 0   & 0       & 0   & 5             & 5   & 0         & 0 & 2420    & 2550 & 0.39      & 1.92  & 57895 & 57835 & 4690 & 3571 \\
$E2_5$    & 1           & 0   & 0       & 0   & 7             & 7   & 0         & 0 & 2450    & 2220 & 0.36      & 1.96  & 57655 & 57655 & 3343 & 971  \\
$E2_6$    & 1           & 0   & 0       & 0   & 6             & 7   & 0         & 0 & 2630    & 2430 & 1.11      & 3.34  & 57775 & 57775 & 3314 & 412  \\
$E2_7$    & 4           & 6   & 0       & 0   & 8             & 7   & 0         & 0 & 3040    & 5490 & 0.66      & 1.36  & 57530 & 57275 & 2999 & 1792 \\
$E2_8$    & 1           & 0   & 0       & 0   & 5             & 5   & 0         & 0 & 2140    & 2140 & 0.29      & 1.23  & 57835 & 57835 & 1176 & 1162 \\
$E2_9$    & 3           & 2   & 0       & 0   & 7             & 7   & 0         & 0 & 2630    & 3260 & 0.38      & 0.87  & 57470 & 57480 & 3972 & 1392 \\
$E2_{10}$ & 3           & 0   & 0       & 0   & 6             & 3   & 0         & 0 & 4000    & 6310 & 0.47      & 1.61  & 57355 & 57404 & 6290 & 4279 \\ \midrule
Avg.          & 2           & 1   & 0       & 0   & 5.7           & 5.4 & 0         & 0 & 2613    & 3235 & 0.50      & 1.99  & 57690 & 57649 & 3850 & 1982 \\ \bottomrule 
\end{tabular}
\end{table}

\begin{table}[htbp]
\footnotesize

\centering
\caption{Solution characteristics at the end of day of execution for 7 flight copies}
\begin{tabular}{l|r|lr|lr|lr|lr|lr|lr|l}
    \FL  \multicolumn{1}{c}{} & \multicolumn{7}{r|}{RCRP } &  \multicolumn{7}{l}{TCRP}
    \ML Instance & \multicolumn{2}{c}{Alt} & \multicolumn{2}{c}{Ret} & \multicolumn{2}{c}{Res} & \multicolumn{2}{c}{Canc} & \multicolumn{2}{c}{Pen} & \multicolumn{2}{c}{Costs} & \multicolumn{2}{c}{Time} \\ \midrule
$E2_1$    & 9           & 14  & 3       & 2   & 11            & 11  & 0         & 0 & 5460    & 7590 & 57788 & 58014 & 673  & 1111 \\
$E2_2$    & 0           & 5   & 0       & 0   & 12            & 9   & 0         & 0 & 2800    & 5750 & 57895 & 57429 & 50   & 203  \\
$E2_3$    & 2           & 5   & 0       & 0   & 10            & 9   & 0         & 0 & 2560    & 3180 & 57775 & 57590 & 108  & 158  \\
$E2_4$    & 3           & 5   & 0       & 0   & 8             & 7   & 0         & 0 & 3380    & 4800 & 57605 & 57470 & 134  & 327  \\
$E2_5$    & 3           & 8   & 0       & 2   & 12            & 11  & 0         & 0 & 4610    & 5680 & 57544 & 57994 & 209  & 1094 \\
$E2_6$    & 2           & 4   & 0       & 0   & 10            & 10  & 0         & 0 & 3630    & 3940 & 57494 & 57550 & 37   & 350  \\
$E2_7$    & 0           & 1   & 0       & 0   & 10            & 9   & 0         & 0 & 3260    & 3480 & 57535 & 57610 & 57   & 166  \\
$E2_8$    & 2           & 3   & 0       & 0   & 10            & 10  & 0         & 0 & 3580    & 3910 & 57765 & 57765 & 79   & 152  \\
$E2_9$    & 0           & 1   & 0       & 0   & 10            & 11  & 0         & 0 & 3580    & 3750 & 57319 & 57355 & 170  & 233  \\
$E2_{10}$ & 2           & 6   & 0       & 0   & 7             & 7   & 0         & 0 & 6120    & 7790 & 57415 & 57175 & 138  & 359  \\ \midrule
Avg.          & 2.3         & 5.2 & 0.3     & 0.4 & 10            & 9.4 & 0         & 0 & 3898    & 4987 & 57614 & 57595 & 166  & 415   \\ \bottomrule   \end{tabular}
\label{tab:}
\end{table}

\begin{table}[htbp]
\footnotesize

\centering
\caption{Solution characteristics at the end of day of execution for 9 flight copies}
\begin{tabular}{l|r|lr|lr|lr|lr|lr|lr|l}
    \FL  \multicolumn{1}{c}{} & \multicolumn{7}{r|}{RCRP } &  \multicolumn{7}{l}{TCRP}
    \ML Instance & \multicolumn{2}{c}{Alt} & \multicolumn{2}{c}{Ret} & \multicolumn{2}{c}{Res} & \multicolumn{2}{c}{Canc} & \multicolumn{2}{c}{Pen} & \multicolumn{2}{c}{Costs} & \multicolumn{2}{c}{Time} \\ \midrule
$E2_1$    & 9           & 13  & 0       & 3   & 12            & 10  & 0         & 0 & 5430    & 7920 & 57619 & 57429 & 1638 & 3632 \\
$E2_2$    & 0           & 4   & 0       & 0   & 12            & 10  & 0         & 0 & 2800    & 4880 & 57895 & 57524 & 209  & 658  \\
$E2_3$    & 2           & 5   & 0       & 0   & 10            & 9   & 0         & 0 & 2720    & 3180 & 57775 & 57590 & 533  & 610  \\
$E2_4$    & 3           & 6   & 0       & 0   & 8             & 7   & 0         & 0 & 3420    & 4360 & 57605 & 57420 & 328  & 940  \\
$E2_5$    & 4           & 7   & 0       & 0   & 11            & 10  & 0         & 0 & 4310    & 5110 & 57655 & 57504 & 1356 & 1376 \\
$E2_6$    & 2           & 4   & 0       & 0   & 10            & 10  & 0         & 0 & 3570    & 3800 & 57494 & 57550 & 116  & 989  \\
$E2_7$    & 0           & 1   & 0       & 0   & 10            & 8   & 0         & 0 & 3140    & 4600 & 57535 & 57419 & 187  & 377  \\
$E2_8$    & 2           & 3   & 0       & 0   & 10            & 11  & 0         & 0 & 3880    & 4280 & 57765 & 57765 & 417  & 813  \\
$E2_9$    & 0           & 1   & 0       & 0   & 10            & 11  & 0         & 0 & 3760    & 3870 & 57374 & 57355 & 359  & 710  \\
$E2_{10}$ & 2           & 7   & 0       & 0   & 7             & 7   & 0         & 0 & 5890    & 7880 & 57665 & 57135 & 935  & 964  \\ \midrule
Avg.      & 2.4         & 5.1 & 0       & 0.3 & 10            & 9.3 & 0         & 0 & 3892    & 4988 & 57638 & 57469 & 608  & 1107 \\ \bottomrule  \end{tabular}
\end{table}

\begin{table}[htbp]
\footnotesize

\centering
\caption{Solution characteristics at the end of day of execution for 11 flight copies}
\begin{tabular}{l|r|lr|lr|lr|lr|lr|lr|l}
    \FL  \multicolumn{1}{c}{} & \multicolumn{7}{r|}{RCRP } &  \multicolumn{7}{l}{TCRP}
    \ML Instance & \multicolumn{2}{c}{Alt} & \multicolumn{2}{c}{Ret} & \multicolumn{2}{c}{Res} & \multicolumn{2}{c}{Canc} & \multicolumn{2}{c}{Pen} & \multicolumn{2}{c}{Costs} & \multicolumn{2}{c}{Time} \\ \midrule
$E2_1$    & 10          & 18  & 0       & 2   & 12            & 10  & 0         & 0 & 5750    & 7840 & 57675 & 57509 & 3606 & 3747 \\
$E2_2$    & 0           & 4   & 0       & 0   & 12            & 10  & 0         & 0 & 2910    & 4480 & 57895 & 57524 & 818  & 2309 \\
$E2_3$    & 3           & 5   & 0       & 0   & 9             & 9   & 0         & 0 & 3320    & 3420 & 57570 & 57590 & 834  & 1453 \\
$E2_4$    & 3           & 5   & 0       & 0   & 8             & 7   & 0         & 0 & 3560    & 4870 & 57965 & 57470 & 3620 & 1650 \\
$E2_5$    & 4           & 7   & 0       & 0   & 11            & 10  & 0         & 0 & 4230    & 5300 & 57655 & 57509 & 3610 & 3479 \\
$E2_6$    & 2           & 5   & 0       & 0   & 10            & 9   & 0         & 0 & 3480    & 4020 & 57494 & 57324 & 357  & 1819 \\
$E2_7$    & 0           & 0   & 0       & 0   & 10            & 8   & 0         & 0 & 3120    & 5760 & 57530 & 57245 & 361  & 1311 \\
$E2_8$    & 2           & 3   & 0       & 0   & 10            & 10  & 0         & 0 & 3750    & 3920 & 57765 & 57765 & 876  & 1612 \\
$E2_9$    & 0           & 1   & 0       & 0   & 10            & 11  & 0         & 0 & 3630    & 4280 & 57374 & 57355 & 1453 & 2797 \\
$E2_{10}$ & 4           & 6   & 0       & 0   & 7             & 7   & 0         & 0 & 13812   & 8010 & 58689 & 57210 & 1072 & 2463 \\ \midrule
Avg.          & 2.8         & 5.4 & 0       & 0.2 & 9.9           & 9.1 & 0         & 0 & 4756.2  & 5190 & 57761 & 57450 & 1661 & 2264 \\ \bottomrule 
   \end{tabular}
\label{tab:comp6}
\end{table}

\section{Conclusions} \label{sec:conclusion}

The increase in airline traffic is forcing the airline companies to deal with disruptions in a more effective way. Airline's capability to respond to these disruptions largely depend on its crew member scheduling. Reserve crew play an important role in disruption management. It is therefore important to design a reliable crew member scheduling by using reserve-crew to mitigate the risks. 

In this paper, we propose a formulation for the reliable crew recovery problem (RCRP), where the goal is to recover a disrupted crew schedule in a cost-efficient manner while remaining capable of recovering from further disruptions. The RCRP explicitly examines the effect of using reserve crew on the resulting robustness of the schedule. We measure the effect on schedule robustness in terms of the so-called expected shortfall. We model the underlying reserve schedules using Markov-chains. We solve the problem using a novel branch-and-price approach. New pairings are created using a variant of a pulse algorithm to solve a Resource Constrained Shortest Path Problem. 

Experiments on real-life data from a medium-sized Dutch carrier show that the RCRP outperforms traditional recovery models in delivering a more stable schedule for the day of execution, which leads to a reduced amount of last-minute crew alterations (and subsequent delays) and even a reduced amount of cancellations due to lack of crew. This is especially important for a carrier operating in a hub-and-spoke network, in which the reserve crew members are located at the main hub from which the aircraft rotations are departing. 

Some limitations exist in the current study that can be considered as the directions of future research. First, we assume that the demand for reserve-crew members is independently distributed among periods. To model the impact of disruption in a more realistic way, the demands in different periods may be correlated to each other. To tackle this challenge, some forecasting models could be used to model the demand correlation. Second, as an initial study of reliable reserve-crew scheduling, we mainly investigate the recovery plan for a single execution day. It may be interesting to consider the plan for multiple flight execution days. Then, a more detailed model is required to track which individual crew members are deadheading. 
{Finally, this study is centered around the development of a reliable scheduling policy at the operational level. Notably, it intentionally confines itself to addressing operational intricacies, without delving into the broader tactical or strategic concerns such as the optimal number of reserve crew members. Consequently, there emerges a significant and compelling avenue for the creation of an integrated planning framework. This framework would inherently possess the capability to seamlessly explore and harmonize both operational nuances and tactical/strategic decisions.}

\bibliographystyle{informs2014trsc}

\bibliography{thesis}

\newpage
\renewcommand{\thesection}{A}
\section{Appendix} \label{sec:appendixA}
\setcounter{section}{0}
\renewcommand{\thesection}{A}
\setcounter{table}{0}
\renewcommand{\thetable}{A\arabic{table}}
\setcounter{figure}{0}
\renewcommand{\thefigure}{A\arabic{figure}}

 	\subsection{Traditional Crew Recovery Problem (TCRP)} \label{sec:TCRPformulation}
\begin{alignat}{2}
\text{minimize } & \sum_{k \in K} \sum_{p\in P^k}c_p^k x_p^k + \sum_{\omega \in \Omega}c^{d}_{\omega}y_{\omega} + \sum_{f\in F}c^{c}_f z_f  ,\\ 
 \text{subject to } & \sum_{k\in K} \sum_{p \in P^k} \sum_{\omega \in \Omega}a^{k}_{\omega p} x_p^k - \sum_{\omega\in\Omega}y_{\omega} + z_f = 1 &\ & \quad \forall f\in \cal{F}, \\
 &  y_{\omega} - \sum_{k \in \cal{K}} \sum_{p \in \cal{P}^k} a_{\omega p}^kx_p^k \leq 0   &\ & \quad \forall \omega \in \Omega \label{lp:deadheadconsistency} \\
 &  \sum_{w \in \Omega_f} y_{\omega} \leq M(1-z_f)   &\ & \quad \forall f \in \cal{F} \label{lp:deadheadconsistency2} \\
&  \sum_{k\in K}\sum_{p\in P^k}(a_{\omega p}^{k}x_p^k - y_\omega) +  \sum_{k\in K}\sum_{p\in P^k}(a_{\omega'p}^{k}x_p^k - y_{\omega'}) \le 1&\ & \quad \forall (\omega, \omega') \in \cal{I},  \\
& \sum_{p \in P^k} x_p^k \le 1 &\ & \quad \forall k\in \cal{K}\backslash \cal{K}^R,  \\
& \sum_{p \in P^k} x_p^k \le 1 &\ & \quad \forall k\in \cal{K}^R,  \\
& \sum_{k\in K^R}\sum_{p \in P^k} x_p^k \le n, &\ &   \\
& x_{p}^k \in \{0,1\} &\ & \quad \forall k\in \cal{K}, \enspace \forall p\in \cal{P}^k,  \\
& y_\omega \ge 0 &\ & \quad \forall \omega \in \Omega,   \\
& \upsilon_{\omega}^+ \ge 0, \upsilon_{\omega}^- \ge 0 &\ & \quad \forall \omega \in \Omega. 
\end{alignat}

\noindent Here, $n$ denotes a limit on the maximum number of reserve crew. 
For a full description of the remaining parameters, we refer the reader to Section \ref{sec:propmodel}.
\end{document}